\documentclass[leqno,12pt]{article}

\usepackage{amstext    }
\usepackage{amsthm    }
\usepackage{a4}
\usepackage[mathscr]{eucal}
\usepackage{mathrsfs}

\usepackage{amsmath}
\usepackage{amssymb}
\usepackage{amscd}

\numberwithin{equation}{section}

\newtheorem{theorem}{Theorem}[section]
\newtheorem{definition}[theorem]{Definition}
\newtheorem{proposition}[theorem]{Proposition}
\newtheorem{corollary}[theorem]{Corollary}
\newtheorem{lemma}[theorem]{Lemma}
\newtheorem{remark}[theorem]{Remark}
\newtheorem{remarks}[theorem]{Remarks}
\newtheorem{example}[theorem]{Example}

\newcommand{\cali}[1]{\mathscr{#1}}

\newcommand{\GL}{{\rm GL}}

\newcommand{\Aut}{{\rm Aut}}

\newcommand{\supp}{{\rm supp}}

\newcommand{\const}{\mathop{\mathrm{const}}}

\newcommand{\dist}{\mathop{\mathrm{dist}}\nolimits}

\newcommand{\loc}{{loc}}
\newcommand{\ddc}{dd^c}
\newcommand{\dc}{d^c}

\newcommand{\dbar}{{\overline\partial}}
\newcommand{\ddbar}{{\partial\overline\partial}}
\renewcommand{\GL}{{\rm GL}}

\newcommand{\DSH}{{\rm DSH}}
\newcommand{\PGL}{{\rm PGL}}

\newcommand{\id}{{\rm id}}

\newcommand{\Kob}{{\rm Kob}}
\newcommand{\adeg}{{\rm deg}_a}

\newcommand{\vol}{{\rm vol}}
\newcommand{\reg}{{\rm reg}}

\newcommand{\Cc}{\cali{C}}
\newcommand{\Dc}{\cali{D}}

\newcommand{\Fc}{\cali{F}}
\newcommand{\Gc}{\cali{G}}

\newcommand{\Jc}{\cali{J}}
\newcommand{\Kc}{\cali{K}}

\newcommand{\Uc}{\cali{U}}

\newcommand{\Tc}{\cali{T}}

\newcommand{\FS}{{\rm FS}}
\newcommand{\B}{\mathbb{B}}
\newcommand{\D}{\mathbb{D}}
\newcommand{\C}{\mathbb{C}}

\newcommand{\Z}{\mathbb{Z}}
\newcommand{\R}{\mathbb{R}}

\renewcommand\P{\mathbb{P}}

\title{Rigidity of Julia sets for H\'enon type maps}

\author{Tien-Cuong Dinh and Nessim Sibony}

\begin{document}

\maketitle

\begin{abstract}
We prove that the Julia set of a H\'enon type automorphism on $\C^2$ is very rigid: it supports a unique positive $\ddc$-closed current of mass 1. 
A similar property holds for the cohomology class of the Green current associated with an automorphism of positive entropy on a compact K\"ahler surface. 
Relations between this phenomenon, several quantitative equidistribution properties and  the theory of value distribution will be discussed. We also survey some rigidity properties of H\'enon type maps on $\C^k$ and of automorphisms of compact K\"ahler manifolds. 
\end{abstract}

\noindent
{\bf Classification AMS 2010:} 37-02, 37F10, 32H30, 32H50, 32U90.


\tableofcontents

\section{Introduction} \label{intro}

The aim of these notes is to explore a rigidity phenomenon for polynomial automorphisms of $\C^k$ and also for holomorphic automorphisms of compact K\"ahler manifolds. This property plays a central role in the dynamical study of these maps and can be seen as a strong ergodicity property in the complex setting. 

Consider first, for simplicity, a polynomial automorphism in $\C^2$ of H\'enon type
$$f(z_1,z_2)=(p(z_1)+az_2,z_1)$$
where $p$ is a one variable polynomial of degree $d\geq 2$ and $a\in \C^*$. Denote by $f^n$ the iterate of order $n$ of $f$ and define 
$$K_+:=\big\{z\in\C^2,\ (f^n(z))_{n\geq 0} \text{ is bounded in } \C^2\big\}.$$

It was shown by J.-E. Forn\ae ss and the second author that $K_+$
supports a unique positive closed $(1,1)$-current of mass 1, the so-called  Green curent $T_+$ \cite{FS1}, see also \cite{DS6}. 
Here we show that indeed $T_+$ is the unique positive $\ddc$-closed $(1,1)$-current of mass 1 with support in  $K_+$. As a direct consequence, we get that if $\tau_n$ are positive $(1,1)$-currents of mass 1 such that $\supp(\tau_n)$ converge to a subset of $K_+$ and $\ddc \tau_n$ converge to 0 then $\tau_n$ converge to $T_+$.
This observation permits to explore the rigidity property and deduce dynamical properties of $f$.

Consider for example any non-constant holomorphic map 
$$\phi:\C\to K_+$$
and the currents of integration on the images of discs of center 0 and of radius $\leq r$ in $\C$. Using a classical idea from Nevanlinna theory, one can average these currents with appropriate weight in order to get a family of currents $\tau_r$ (Nevanlinna's currents) satisfying the above properties, see Section \ref{section_rigid} for details. They necessarily converge to $T_+$ as $r\to\infty$. One can replace $\C$ with a parabolic Riemann surface or a disc on which $\phi$ satisfies an appropriate growth condition. 

What is striking in the previous result is the claim that  the currents $\tau_r$ not only converge, but that the limit is somehow independent of $\phi$. The result applies for example when $\phi$ parametrizes a stable manifold associated with $f$ and gives us a rough information about the behavior of this stable manifold.

We then extend the result to holomorphic automorphisms of compact K\"ahler surfaces. Let $f:X\to X$ be an automorphism of positive entropy of a compact K\"ahler surface $X$. We will show that the Green $(1,1)$-current $T_+$ associated with $f$ is the unique positive $\ddc$-closed $(1,1)$-current in its cohomology class $\{T_+\}$. This current is positive and closed, see \cite{Cantat1,DS1,DS2}.  
 In this case, if a holomorphic map $\phi:\C\to X$ satisfies $\phi^*(T_+)=0$, then Nevanlinna's currents $\tau_r$ associated with $\phi$ also converge to $T_+$ as $r\to\infty$. 
 
The condition $\phi^*(T)=0$ is realized in particular when a subsequence of $(f^{n})_{n\geq 0}$ is locally equicontinuous on the image $\phi(\C)$. In the case where $T_+$ can be expressed on an open set as an average of currents of integration on disjoint Riemann surfaces, then
$\phi^*(T)=0$ means that the image of $\phi$ is along these Riemann surfaces. 
So if we consider that the family of maps $\phi:\C\to X$ with $\phi^*(T)=0$ is a "$T$-lamination" we get a unique ergodicity result for that $T$-lamination. 
A quantitative version of the rigidity of the class $\{T_+\}$ is given in Theorem \ref{th_explicit_rigidity}. It says that the distance between a positive closed $(1,1)$-current $S$ and $T_+$ is bounded by an explicit function of the distance between the classes $\{S\}$ and  $\{T_+\}$. 

This point of view provides a strong analogy between equidistribution properties for stable manifolds in discrete holomorphic dynamics and ergodic properties of foliations by Riemann surfaces as developed in \cite{DNS1, FS2,FS3}.
We however do not discuss here the theory of foliations. Though the dictionary between the theory of discrete holomorphic dynamical systems in several variables and the theory of foliations deserves to be explored further.

We also consider in this paper several equidistribution properties towards the Green currents with precise control of the convergence speed. A large class of polynomial automorphisms on $\C^k$ and automorphisms of higher dimensional compact K\"ahler manifolds are also studied. 

The plan of the article is as follows. In Section \ref{section_psh}, we recall basic results in pluripotential theory that will be used later. Further results for currents in compact K\"ahler manifolds are given in Section \ref{section_Kahler}. The notions of rigid set and rigid cohomology class are introduced in Section \ref{section_rigid}. H\'enon maps are considered in Sections \ref{section_Henon} and \ref{section_Green} while automorphisms of compact K\"ahler surfaces are treated in Section \ref{section_surface} and the higher dimensional dynamics in Section \ref{section_higher_dim}. 
It is possible to start reading from Section \ref{section_Henon} and come back to the technical tools developed in Sections \ref{section_psh}, \ref{section_Kahler} and \ref{section_rigid}, in particular to the important estimate given in Corollary \ref{cor_dsh_exp}. 

Finally, we refer to 
Bedford-Lyubich-Smillie, de Th\'elin and Dinh-Sibony \cite{BLS1, BLS2, BS1, deThelin2,DS9} for the ergodic properties of the measures of maximal entropy and the distribution of periodic points in the case of H\'enon type maps. 
The Green currents and the measure of maximal entropy were introduced by the second author of the present paper as noticed in \cite[p.78]{BLS1}. 
For  automorphisms of compact K\"ahler manifolds see Cantat, de Th\'elin-Dinh, Dinh-Sibony \cite{Cantat1,dD,DS1,DS2,DS9}, for
the semi-local setting of horizontal-like maps see Dinh-Nguyen-Sibony, Dujardin \cite{DNS2,DS4,Du}, for non-injective holomorphic maps see the survey \cite{DS7} and the references therein.

\medskip
\noindent
{\bf Acknowledgment.} We would like to thank the referee for his remarks which permitted us to improve the presentation of this article.

\section{Positive closed currents and p.s.h. functions} \label{section_psh}

In this section, we recall briefly the notions of positive closed currents and of plurisubharmonic (p.s.h. for short) functions on a complex manifold. We refer to \cite{Demailly3,DS7} for an account of that theory. The readers who are familiar with  pluripotential theory may skip this section. 

\medskip

\noindent
{\bf $\bullet$ Differential forms on complex manifolds.} 
Let $X$ be a complex manifold of dimension $k$, e.g. $\C^k$ or an open subset of $\C^k$. Let $\varphi$ be a
differential $l$-form on $X$. 
In local holomorphic
coordinates $z=(z_1,\ldots,z_k)$, it can be written as
$$\varphi(z)=\sum_{|I|+|J|=l}\varphi_{IJ} dz_I\wedge d\overline z_J,$$
where $\varphi_{IJ}$ are complex-valued functions,
$dz_I:=dz_{i_1}\wedge\ldots\wedge dz_{i_p}$ if
$I=(i_1,\ldots,i_p)$, and  $d\overline
z_J:=d\overline z_{j_1}\wedge\ldots\wedge d\overline z_{j_q}$ if $J=(j_1,\ldots,j_q)$. 

The {\it conjugate} of
$\varphi$ is defined by
$$\overline\varphi(z):=\sum_{|I|+|J|=l}\overline \varphi_{IJ} d\overline z_I\wedge d z_J.$$
The form $\varphi$ is {\it real} if and only if $\varphi=\overline\varphi$. 

We say that $\varphi$ is a form of {\it of bidegree} $(p,q)$ if
$\varphi_{IJ}=0$ when $(|I|,|J|)\not=(p,q)$. The bidegree does not
depend on the choice of local holomorphic coordinates.
Let $T_X^\C$ denote the complexification of the tangent bundle of
$X$. The complex structure on $X$ induces a linear endomorphism $\Jc$ on
the fibers of $T_X^\C$ such that $\Jc^2=-\id$. This endomorphism
induces a decomposition of 
$T_X^\C$ into the direct sum of two proper sub-bundles of dimension $k$: {\it the
  holomorphic part} $T_X^{1,0}$ associated with the eigenvalue $\sqrt{-1}$ of
$\Jc$, and the {\it anti-holomorphic part}
$T_X^{0,1}$ associated with the eigenvalue $-\sqrt{-1}$.
Let $\Omega_X^{1,0}$ and $\Omega_X^{0,1}$ denote the dual bundles
of $T_X^{1,0}$  and  $T_X^{0,1}$. We can consider $(p,q)$-forms as sections of
the vector bundle $\bigwedge^p\Omega^{1,0}\otimes \bigwedge^q\Omega^{0,1}$.

If $\varphi$ is a
$(p,q)$-form then the differential $d\varphi$ is the sum of a $(p+1,q)$-form 
and a $(p,q+1)$-form. We  denote by $\partial \varphi$ the part of
bidegree $(p+1,q)$ and  by  $\dbar \varphi$ the  the part of
bidegree $(p,q+1)$. The operators $\partial$ and $\dbar$ extend
linearly to arbitrary forms $\varphi$. 
The operator $d$ is real, i.e. it
sends real forms to real forms but
$\partial$ and $\dbar$ are not
real. The identity $d\circ d=0$ implies that
$\partial\circ\partial=0$, $\dbar\circ\dbar=0$ and $\partial\dbar+\dbar\partial=0$. 

Define $\dc:={\sqrt{-1}\over2\pi}(\dbar-\partial)$. This operator is real
and satisfies $\ddc = {\sqrt{-1}\over \pi}\ddbar$. 
Note that the above operators 
commute with the pull-back by holomorphic maps. More precisely, if
$\tau:X_1\rightarrow X_2$ is a holomorphic map between complex
manifolds and $\varphi$ is a form on $X_2$ then
$d\tau^*(\varphi)=\tau^*(d\varphi)$, $\ddc \tau^*(\varphi)=\tau^*(\ddc\varphi)$, etc.
Recall that the form $\varphi$ is {\it closed}
(resp. $\partial$-closed, $\dbar$-closed, $\ddc$-closed) if $d\varphi$
(resp. $\partial\varphi$, $\dbar\varphi$, $\ddc\varphi$) vanishes. The
form $\varphi$ is {\it exact} (resp. $\partial$-exact, $\dbar$-exact,
$\ddc$-exact) if it is equal to the differential
$d\psi$ (resp. $\partial \psi$, $\dbar\psi$, $\ddc\psi$) of a form $\psi$. Clearly, exact forms are closed.

A smooth $(1,1)$-form $\omega$ on $X$ is {\it Hermitian} if it
can be written in local coordinates as 
$$\omega(z)=\sqrt{-1}\sum_{1\leq i,j\leq k} \alpha_{ij}(z) dz_i\wedge d\overline z_j,$$ 
where $\alpha_{ij}$ are smooth functions such that the matrix
$(\alpha_{ij})$ is Hermitian. We consider a form $\omega$ such that
the matrix $(\alpha_{ij})$ is
positive definite at every point. It is
strictly positive in the sense that we will introduce later. The form $\omega$ induces a so-called {\it Hermitian metric} on $X$ as follows. 

The form $\omega$ is always real
and induces a norm on the tangent spaces of $X$. So it defines a
Riemannian metric on $X$. For example, the Euclidean metric on $\C^k$ is associated with the standard Hermitian 
form
$$\beta:=\sqrt{-1}\sum_{1\leq i\leq k} dz_i\wedge d\overline z_i.$$ 
In general, for each point $a\in X$,  we can choose local coordinates $z$ near $a$ such that $z=0$ at $a$ and 
$$\omega=\sqrt{-1}\sum_{1\leq i\leq k} dz_i\wedge d\overline z_i \quad \text{at} \quad a.$$ 
So the Riemannian metric associated with $\omega$ coincides at $a$ with the Euclidean metric on the above chart. 
It is easy to construct Hermitian metrics on $X$ using local coordinates and a partition of unity. 
From now on, we assume that $X$ is endowed with a fixed Hermitian metric $\omega$.

The following result is due to Wirtinger. The remarkable fact is that in order to compute the volume of an analytic set, we have to integrate a form which is independent of the analytic set. 

\begin{theorem}[Wirtinger] \label{th_wirtinger}
Let $Y$ be an analytic set of pure dimension
  $p$ in a Hermitian manifold $(X,\omega)$. Then the $2p$-dimensional
  volume of $Y$ in a Borel set $K$ is equal to
$$\vol(Y\cap K)={1\over p!}\int_{\reg(Y)\cap K} \omega^p.$$
Here, the volume is with respect to the Riemannian metric induced by
$\omega$ and $\reg(Y)$ denotes the set of regular points in $Y$. 
\end{theorem}

\smallskip
\noindent
{\bf $\bullet$ Currents on a complex manifold.}
We now  introduce positive forms and positive
currents  
on complex manifolds. The concept of positivity is due to Lelong and Oka. The theory has many
applications in complex geometry, algebraic geometry and dynamics, see \cite{Demailly3, DS6}.

Let $(X,\omega)$ be a Hermitian manifold of dimension $k$. 
Recall that {\it a current} $S$ on $X$, of degree $l$ and of dimension $2k-l$, 
is a continuous linear form on the
space $\Dc^{2k-l}(X)$ of smooth $(2k-l)$-forms with compact
support in $X$. Its value on a $(2k-l)$-form $\varphi\in\Dc^{2k-l}(X)$ is
denoted by $S(\varphi)$ or more frequently by $\langle S,\varphi\rangle$.
On a chart, $S$ corresponds to a continuous linear form acting on
the coefficients of $\varphi$. So it can be represented as an $l$-form with
distribution coefficients. 

A sequence $(S_n)$ of $l$-currents
{\it converges} to an $l$-current $S$ if for every $\varphi\in\Dc^{2k-l}(X)$,
$\langle S_n,\varphi\rangle$ converge to $\langle S,\varphi\rangle$. 
{\it The conjugate} of $S$ is the $l$-current $\overline S$ defined by 
$$\langle\overline S,\varphi\rangle:=\overline{\langle
  S,\overline\varphi\rangle},$$
for $\varphi\in\Dc^{2k-l}(X)$. The current $S$ is {\it real} if and only if
$\overline S=S$.
{\it The support} of $S$ is the smallest closed subset $\supp(S)$ of $X$
such that $\langle S,\varphi\rangle=0$ when $\varphi$ has compact support in
$X\setminus \supp(S)$. The current $S$ extends continuously to the
space of smooth forms $\varphi$ such that $\supp(\varphi)\cap\supp(S)$
is compact in $X$. 

If $\alpha$ is a smooth $s$-form on $X$ with $s\leq 2k-l$, we define the $(l+s)$-current $S\wedge\alpha$ by 
$$\langle S\wedge \alpha,\varphi\rangle:=\langle S,\alpha\wedge\varphi\rangle$$
for every form $\varphi\in \Dc^{2k-l-s}(X)$. Define also $\alpha\wedge S:=(-1)^{ls}S\wedge\alpha$. 

If $X'$ is a complex manifold of dimension $k'$
with $2k'\geq 2k-l$, and if $\tau:X\rightarrow X'$ is a holomorphic map
which is proper on the support of $S$, we can define {\it the
  push-forward} $\tau_*(S)$ of $S$ by $\tau$. The current $\tau_*(S)$ has the
same dimension than $S$, i.e. of degree $2k'-2k+l$ and is supported
on $\tau(\supp(S))$. It satisfies
$$\langle \tau_*(S),\varphi\rangle:=\langle S,\tau^*(\varphi)\rangle$$
for $\varphi\in \Dc^{2k-l}(X')$. 

If $X'$ is a complex manifold of
dimension $k'\geq k$ and if $\tau:X'\rightarrow X$ is a submersion, we
can define {\it the pull-back} $\tau^*(S)$ of $S$ by $\tau$. This is an $l$-current
supported on $\tau^{-1}(\supp(S))$, it satisfies
$$\langle \tau^*(S),\varphi\rangle:=\langle S,\tau_*(\varphi)\rangle$$
for $\varphi\in \Dc^{2k'-l}(X')$. Indeed, since $\tau$ is a
submersion, the current $\tau_*(\varphi)$ is in fact a smooth form
with compact support in $X$; it is given by an integral of $\varphi$
on the fibers of $\tau$.

Any smooth differential $l$-form $\psi$ on $X$ can be considered as a current: it
defines the continuous linear form $\varphi\mapsto
\int_X\psi\wedge\varphi$ on $\varphi\in\Dc^{2k-l}(X)$. So
currents extend the notion of differential forms. 
The operators $d,\partial,\dbar$ on differential forms extend to
currents. For example, we have that $dS$ is an $(l+1)$-current defined
by
$$\langle dS,\varphi\rangle:=(-1)^{l+1}\langle S, d\varphi\rangle$$
for $\varphi\in\Dc^{2k-l-1}(X)$. One easily check that when $S$ is
a smooth form, the above identity is a consequence of the Stokes' formula.

We say that $S$ is of 
{\it bidegree} $(p,q)$ and {\it of bidimension} $(k-p,k-q)$ if it
vanishes on forms of bidegree $(r,s)\not=(k-p,k-q)$. The conjugate of
a $(p,q)$-current is of bidegree $(q,p)$. So, if such a current is real,
we necessarily have $p=q$. Note that the push-forward and the
pull-back by holomorphic maps commute with the above operators. They
preserve real currents; the push-forward preserves the bidimension and
the pull-back preserves the bidegree.

\medskip
\noindent
{\bf $\bullet$ Positive forms and positive currents.}
There exist three notions of positivity which coincide for the bidegrees
$(0,0)$, $(1,1)$, $(k-1,k-1)$ and $(k,k)$. Here, we only use two of
them. They are dual to each other. 

\begin{definition} \rm
A $(p,p)$-form $\varphi$ is {\it (strongly) positive} if at each
point, it is equal to a combination with positive coefficients of forms
of type
$$(\sqrt{-1}\alpha_1\wedge \overline\alpha_1)\wedge \ldots \wedge
(\sqrt{-1}\alpha_p\wedge \overline\alpha_p),$$
where $\alpha_i$ are $(1,0)$-forms. 
\end{definition}

Any $(p,p)$-form can be written as a finite combination of positive $(p,p)$-forms.
For example, in local coordinates $z$, a $(1,1)$-form $\omega$ is written as
$$\omega=\sum_{i,j=1}^k \alpha_{ij} \sqrt{-1}dz_i\wedge d\overline z_j,$$
where $\alpha_{ij}$ are functions. This form is positive if and only
if the matrix $(\alpha_{ij})$ is positive semi-definite at every
point. In local coordinates $z$, the $(1,1)$-form $\ddc\|z\|^2$ is
positive. One can write $dz_1\wedge d\overline z_2$ as a combination
of $dz_1\wedge d\overline z_1$, $dz_2\wedge d\overline z_2$, 
$d(z_1\pm z_2)\wedge d\overline{(z_1\pm z_2)}$ and 
$d(z_1\pm \sqrt{-1}z_2)\wedge d\overline{(z_1\pm \sqrt{-1}z_2)}$.  Hence, we see that positive
forms generate the space of $(p,p)$-forms.

\begin{definition} \rm
Let $S$ be a $(p,p)$-current on $X$. We say that $S$ is {\it
  weakly positive} if  $S\wedge \varphi$ is a positive measure for
every smooth positive $(k-p,k-p)$-form $\varphi$,
and that $S$ is {\it positive}
if $S\wedge \varphi$ is a positive measure for
every smooth weakly positive $(k-p,k-p)$-form $\varphi$. 
\end{definition}

Positivity implies weak positivity. These properties 
 are preserved
under pull-back by holomorphic submersions and push-forward by 
proper holomorphic maps.
Positive and weakly positive forms or
currents are real. One can consider positive and weakly positive
$(p,p)$-forms as sections of some bundles of strictly convex closed
cones  in
the real part of the vector bundle $\bigwedge^p\Omega^{1,0}\otimes
\bigwedge^p\Omega^{0,1}$. 

The wedge-product of a positive current with a positive form is positive.
The wedge-product of a weakly positive current with a positive form is
weakly positive. Wedge-products of weakly positive forms or currents are not
always weakly positive. For real $(p,p)$-currents or
forms $S$, $S'$, we will write $S\geq S'$ and $S'\leq S$ if $S-S'$ is
positive. 

\begin{definition} \rm
A $(p,p)$-current or form $S$ is {\it strictly positive} if in local
coordinates $z$, there is a constant $\epsilon>0$ such that $S\geq
\epsilon (\ddc \|z\|^2)^p$. 
\end{definition}

Equivalently, $S$ is strictly positive if
we have locally $S\geq \epsilon\omega^p$ with $\epsilon>0$.

\begin{example} \rm \label{example_current_lelong}
Let $Y$ be an analytic set of pure codimension $p$ of $X$. 
Using the local description of $Y$ near a singularity 
\cite{Gu} and 
Wirtinger's theorem \ref{th_wirtinger}, one can prove 
that the $2(k-p)$-dimensional volume of $Y$
is locally finite in $X$. This allows to define the following
$(p,p)$-current $[Y]$ by
$$\langle [Y],\varphi\rangle:=\int_{\reg(Y)}\varphi$$
for $\varphi$ in $\Dc^{k-p,k-p}(X)$, the space of smooth
$(k-p,k-p)$-forms with compact support in $X$. Here $\reg(Y)$ denotes the smooth points of $Y$.  Lelong proved that this current is
positive and closed \cite{Demailly3, Lelong}. 
\end{example}

If $S$ is a (weakly) positive $(p,p)$-current, it is of order
0, i.e. it extends continuously to the space of continuous forms with
compact support in $X$. In other words, on a chart of $X$, the current
$S$ corresponds to a differential form with measure coefficients. 

\begin{definition} \rm
The {\it mass} of a positive $(p,p)$-current $S$ on a Borel set $K$ is defined by 
$$\|S\|_K:=\int_K S\wedge\omega^{k-p}.$$
\end{definition}

When $K$ is contained in a fixed compact subset of $X$, we obtain an equivalent
norm if we change the Hermitian metric on $X$. This is a consequence
of an above-mentioned property, which says that $S$ takes values in
strictly convex closed cones. Note that the previous mass-norm is
just defined by an integral, which is easier to compute or
to estimate than the usual mass for currents on real
manifolds. For the current $[Y]$ in Example \ref{example_current_lelong}, by Wirtinger's theorem, the
mass on $K$ is equal to $(k-p)!$ times the volume of $Y\cap K$ with
respect to the considered Hermitian metric.

Positivity implies an important compactness property. As for positive
measures, any family of
positive $(p,p)$-currents with locally uniformly bounded mass, is
relatively compact in the cone of positive $(p,p)$-currents. 
We will need the following classical result.

\begin{theorem} \label{th_ext_current}
Let $E$ be a closed subset of a complex manifold $X$ of dimension $k$. Let $T$ be a positive closed $(p,p)$-current on $X\setminus E$. Assume that the Hausdorff $2(k-p)$-dimensional measure of $E$ vanishes. Then $T$ has finite mass on compact subsets of  $X$ and its extension by $0$ through $E$ is a positive closed $(p,p)$-current on $X$. 
\end{theorem}

\medskip

\noindent
{\bf $\bullet$ Plurisubharmonic functions.}
Calculus on currents is often delicate. However, the theory is
well developped for positive closed $(1,1)$-currents thanks to
plurisubharmonic functions. Note that positive
closed $(1,1)$-currents correspond to hypersurfaces (analytic sets of
pure codimension 1) in complex geometry and working with
$(p,p)$-currents, as with higher codimension analytic sets, is more
difficult.

\begin{definition}\rm 
An upper semi-continuous function $u:X\rightarrow \R\cup\{-\infty\}$,
not identically $-\infty$ on any irreducible component of $X$, is {\it
  plurisubharmonic} (p.s.h. for short) if it is subharmonic or
identically $-\infty$ on
any holomorphic  disc in $X$. A function $v$ is {\it
  pluriharmonic} if $v$ and $-v$ are p.s.h.
\end{definition}

Note that p.s.h. functions are defined at
every point. 
The semi-continuity implies that p.s.h. functions are locally bounded from above.
Pluriharmonic functions
are locally real parts of holomorphic functions, in particular, they
are real analytic. 

Recall that a {\it holomorphic disc} in $X$ is a holomorphic map
$\tau:\Delta\rightarrow X$ where $\Delta$ is the unit disc in
$\C$. One often identifies this holomorphic disc with its image
$\tau(\Delta)$. If $u$ is p.s.h., then $u\circ\tau$ is
subharmonic or identically $-\infty$ on $\Delta$. 

As for subharmonic
functions on $\R^n$, we have the submean inequality: {\it in local
  holomorphic coordinates, the value at $a$ of a p.s.h. function is smaller or
  equal to the average of the function on a sphere centered at
  $a$}. Indeed, this average increases with the radius of the
sphere. 
The submean inequality implies that p.s.h. functions  satisfy the maximum principle: {\it if a
  p.s.h. function on a connected manifold $X$ has a maximum, it is constant}. 
  It also implies that p.s.h. functions are locally integrable. We have the following general properties.
 
\begin{theorem} \label{th_hormander}
P.s.h. functions on $X$ are in $L^p_\loc(X)$ and they form a closed convex cone of $L^p_\loc(X)$  for every $1\leq p<\infty$. 
Let $\Fc$ be a family of p.s.h. functions on $X$ which
  is bounded in $L^1_\loc(X)$.
Then $\Fc$ is relatively compact in $L^p_\loc(X)$ for every $1\leq p<\infty$. Moreover, for every compact subset $K$ of $X$, there are constants $\alpha>0$ and $A>0$ such that
$$\|e^{\alpha |u|}\|_{L^1(K)}\leq A \quad \text{for} \quad u\in\Fc.$$
\end{theorem}

P.s.h. functions are in general unbounded. However, the last
estimate shows that such functions are moderately  unbounded.   
The following propositions are useful in constructing 
p.s.h. functions. 

\begin{proposition} \label{prop_psh_composition}
Let $\chi:\R^n\to\R$ be 
a function which is
convex in all variables and increasing in each variable. Let
$u_1,\ldots,u_n$ be p.s.h. functions on $X$. Then
$\chi(u_1,\ldots,u_n)$ extends through the set $\{u_1=-\infty\} \cup \cdots \cup\{u_n=-\infty\}$ to a p.s.h. function on $X$.  In particular, the function
$\max(u_1,\ldots,u_n)$ is p.s.h. on $X$.
\end{proposition}

\begin{proposition}
Let $E$ be an analytic subset of codimension at least $2$ of  $X$. If $u$ is a p.s.h. function on
$X\setminus E$, then the extension
of $u$ to $X$ given by
$$u(z):=\limsup_{w\rightarrow z\atop w\in X\setminus E} u(w) \quad
\mbox{for}\quad z\in E,$$ 
is a p.s.h. function.  
\end{proposition}

The following result relates p.s.h. functions with positive closed currents.

\begin{proposition}
If $u$ is a p.s.h. function on $X$ then  $\ddc u$ is a positive closed $(1,1)$-current.
Conversely, any positive closed $(1,1)$-current can be locally written
as $\ddc u$ where $u$ is a (local) p.s.h. function. In particular, a smooth function $u$ is p.s.h. if and only if $\ddc u$ is a positive $(1,1)$-form.
\end{proposition}

\begin{definition} \rm
If $S$ is a positive closed $(1,1)$-current on $X$, we call {\it potential} of $S$ any p.s.h. function $u$ such that 
$\ddc u=S$. A p.s.h. function $u$ is called {\it strictly p.s.h.} if the current $S=\ddc
u$ is strictly positive. 
\end{definition}

The above proposition shows that $S$ always admits local potentials.  Two local potentials
differ by a pluriharmonic function. So there is
almost a correspondence between positive closed $(1,1)$-currents and
p.s.h. functions. 
Since pluriharmonic functions are smooth, singularities of positive closed $(1,1)$-currents can be understood via their local potentials.

\begin{example}\rm
Let $f$ be a holomorphic function on $X$ not identically 0 on any
component of $X$. Then,
$\log|f|$ is a p.s.h. function and we have 
$$\ddc \log|f|=
\sum n_i [Z_i],$$ 
where $Z_i$ are irreducible components of the
hypersurface $\{f=0\}$ and
$n_i$ their multiplicities. The last equation is called {\it
  Poincar{\'e}-Lelong equation}. 
  Locally, the ideal of holomorphic
functions vanishing on $Z_i$ is generated by a holomorphic function
$g_i$ and $f$ is equal to the product of $\prod g_i^{n_i}$ with a
non-vanishing holomorphic function. 

In some sense, this is the class of 
the most singular p.s.h. functions. If $X$ is a ball in $\C^k$, 
the convex set generated by such functions is dense in the cone of p.s.h. functions
\cite{Ho, Gu} for the $L^1_\loc$ topology. 
\end{example}

\begin{example} \rm
If $u$ is a p.s.h. function on $X$ and $\tau:X'\to X$ is a holomorphic map, then $u\circ\tau$ is either identically $-\infty$ or a p.s.h. function on each component of $X'$. Since the function $\log\|z\|$ is p.s.h. on $\C^n$, we deduce that   
$\log(|f_1|^2+\cdots+|f_n|^2)$ is p.s.h. on $X$ if
$f_1,\ldots, f_n$ are
holomorphic functions on $X$, not all identically 0 on a component of $X$.  
\end{example}

The following result is useful in the calculus with p.s.h. functions and positive closed $(1,1)$-currents.

\begin{proposition} \label{prop_decreasing_psh}
If $(u_n)$ is a decreasing sequence of p.s.h. functions on $X$, it converges pointwise either to $-\infty$ on at least one component of $X$ or to a p.s.h. function on $X$.  Moreover, every p.s.h. function is locally the limit of a decreasing sequence of smooth p.s.h. functions.
\end{proposition}

\smallskip
\noindent
{\bf $\bullet$ Intersection of currents and slicing.} Let $T$ be a positive closed $(p,p)$-current, $0\leq p\leq k-1$, and $S$ a positive closed $(1,1)$-current on $X$. We will define their wedge-product (intersection) $S\wedge T$. Let $u$ be a local potential of $S$ on an open set $U$  and assume that $u$ is locally integrable on $U$ with respect to the {\it trace measure} $T\wedge\omega^{k-p}$ of $T$. Then the product $uT$ defines a $(p,p)$-current on $U$. We define $S\wedge T$ on $U$ by
$$S\wedge T:=\ddc(uT).$$

It is not difficult to check that if $v$ is a pluriharmonic function then $\ddc(vT)=0$. So the above definition does not depend on the choice of $u$ and then gives a $(p+1,p+1)$-current on $X$ that we denote by $S\wedge T$. 
By definition, this current is locally exact, so it is closed. When $u$ is smooth, this wedge-product is equal to the wedge-product of the positive form $\ddc u$ with $S$ and we see that $S\wedge T$ is positive. The property extends to the general case because we can approximate $u$ by a decreasing sequence of smooth p.s.h. functions, see Proposition \ref{prop_decreasing_psh}. 
Observe also that $\supp(S\wedge T)\subset \supp(S)\cap \supp(T)$. 
We have the following result.

\begin{proposition}
If $(u_n)$ is a sequence of p.s.h. functions decreasing to a p.s.h. function which is locally integrable with respect to the trace measure of $T$, then $\ddc u_n\wedge T$ converge weakly to $\ddc u\wedge T$. If $S$ is a current with local continuous potentials, then $S\wedge T$ depends continuously on $T$.  
\end{proposition}

So we can define by induction the wedge-product of $T$ with several positive closed $(1,1)$-currents. For example, if $S_1,\ldots,S_q$, $q\leq k-p$, are positive closed $(1,1)$-currents with continuous local potentials, then the wedge-product 
$$S_1\wedge\ldots\wedge S_q\wedge T$$
is a positive closed $(p+q,p+q)$-current which depends continuously on $T$ and is symmetric with respect to the currents $S_i$. 

We now consider a special case of the slicing theory that we will use later. Consider a holomorphic map $\pi:X\to Y$ where $Y$ is a complex manifold of dimension $l<k$. 
For simplicity, assume that $\pi$ is a submersion but we can also treat in the same way the case where all fibers of $\pi$ have dimension $k-l$.  
So for every $y\in Y$ the fiber $\pi^{-1}(y)$ is a submanifold of dimension $k-l$ of $X$. Let  $S_1,\ldots,S_q$, $q\leq k-l$, be positive closed $(1,1)$-currents with continuous local potentials on $X$. Define $S:=S_1\wedge \ldots\wedge S_q$ and the slices $\langle S|\pi|y\rangle$ by 
$$\langle S|\pi|y\rangle := S\wedge [\pi^{-1}(y)]=S_1\wedge \ldots\wedge S_q \wedge [\pi^{-1}(y)].$$
We deduce from the above discussion that $\langle S|\pi|y\rangle$ depends continuously on $y$ and its support is contained in $\supp(S)\cap \pi^{-1}(y)$. The continuity of the potentials may be weakened but then we loose the continuous dependence of the slice in $y$.

Observe that different currents may have the same slices. For example, if $\pi:\C^2\to \C$ is the natural projection on the first factor then the slices of the positive closed $(1,1)$-form $T=i dz_1\wedge d\overline z_1$ on $\C^2$ are zero but $T$ does not vanish. We have the following result.

\begin{theorem} \label{th_slice}
Let $\pi:X\to Y$ be a holomorphic submersion as above. Let $T$ and $S$ be positive closed $(1,1)$-currents on $X$ such that $T-S=\ddc u$ for some locally integrable function $u$ on $X$. Assume  that no connected component of any fiber of $\pi$ is contained in $\supp(T)\cup\supp(S)$. Assume also that $\langle T|\pi|y\rangle =\langle S|\pi|y\rangle$ for almost every $y\in Y$. Then $T=S$. 
\end{theorem}
\proof
Let $U$ be the largest open subset of $X$ on which  $u$ is equal to a pluriharmonic function. It is enough to show that $U=X$. Considering a neighbourhood of a point in the boundary of $U$ allows to reduce the problem to the following local setting. We can assume that $\pi:X\to Y$ is the canonical projection from a polydisc $\D^k$ to the first $l$ factors $\D^l$.  In this case, $U$ is an open set in $\D^k$ which intersects all fibers of $\pi$. The function $u$ is pluriharmonic on $U$ and its restriction to almost every fiber of $\pi$ is pluriharmonic. It follows from a version of Levi's extension theorem  that $u$ is equal almost everywhere to a pluriharmonic function. We conclude that $U=X$. This completes the proof of the theorem. 
\endproof

Note that the existence of $u$ in the hypotheses of the above theorem is satisfied when $X$ is a compact K\"ahler manifold and $T$ and $S$ are in the same cohomology class. The positivity of $T$ and $S$ can be replaced by the weaker condition that these currents are of order 0 or flat. In this case, Federer's theory allows to define the slices with respect to $\pi$. 



Recall that a map $\phi$, which associates  $y\in Y$ with a closed subset $F_y$ of $X$, is called {\it upper semi-continuous} (with respect to the Hausdorff topology for closed subsets of $X$)
if for any $y_0\in Y$,
any neighbourhood $V$ of $F_{y_0}$ and for a compact subset $K$ of $X$, we have $F_y\cap K\subset V$ for $y$ close enough to $y_0$. We also say that $\phi$ is {\it lower semi-continuous}  if given $\epsilon>0$ and $y_0\in Y$ we have $\dist(x,F_y)<\epsilon$ for every $x\in F_{y_0}$ and for $y$ close enough to $y_0$. 
The map $\phi$ is {\it continuous} if it is upper and lower semi-continuous.
We have the following elementary lemma.

\begin{lemma}
Let $S=S_1\wedge\ldots\wedge S_q$ be as above. 
The map $y\mapsto \supp(S)\cap \pi^{-1}(y)$ is upper semi-continuous and the map $y\mapsto \supp \langle S|\pi|y\rangle $ is lower semi-continuous with respect to the Hausdorff topology for closed subsets of $X$. 
\end{lemma}
\proof
It is not difficult to see that the first assertion is true if we replace $\supp(S)$ by any closed subset of $X$ and the second assertion holds for any continuous family of currents.
\endproof

The support of $\langle S|\pi|y\rangle$ is contained in $\supp(S)\cap\pi^{-1}(y)$. We will be concerned with the case where these sets are equal.  

\begin{definition} \rm
We call {\it bifurcation locus of $S$ with respect to $\pi$} the closure of the following set
$$\Big\{y\in Y, \quad  \supp(S)\cap \pi^{-1}(y)\not =\supp \langle S|\pi|y\rangle \Big\}.$$
\end{definition}

\smallskip
\noindent
{\bf $\bullet$ Intersection with a positive $\ddc$-closed current.} 
Assume now that $T$ is a positive $\ddc$-closed current, $0\leq p\leq k-1$, and $S$ a positive closed $(1,1)$-current on a complex manifold $X$. Assume that the local potentials of $S$ are continuous on a neighbourhood of the support of $T$. 

We want to define the wedge-product $T\wedge S$ as a current. 
We do this locally and need to check that the definition does not depend on local coordinates. It follows that the definition extends to all complex manifolds.
So we can assume that $X$ is the unit ball $\B$ in $\C^k$ and $S=\ddc u$ where $u$ is a bounded and continuous p.s.h. function on $\B$. We want to define $T\wedge S$ in a neighbourhood of 0. Replacing $u$ by $\max(\log\|z\|,u-c)$ for some constant $c>0$ large enough, we only modify $u$ outside a neighbourhood of 0. This allows to assume that $u=\log\|z\|$ on $\{z\in\C^k,\ r<\|z\| <1\}$ for some constant $0<r<1$. 

Let $\Fc$ denote the set of continuous p.s.h. functions on $\B$ satisfying the last property for a fixed $r$. By a standard method of regularization, functions in $\Fc$ can be uniformly approximated by smooth ones in the same class. In particular, there are $u_n$ smooth in $\Fc$ which converge uniformly to $u$.
 In order to define $T\wedge S$, we only have to show that $T\wedge S_n$, where $S_n:=\ddc u_n$, converge to a current which does not depend on the choice of $u_n$. We then define $T\wedge S$ as equal to the above limit. 
 
 The details and the proof of Theorem \ref{th_wedge_ddc} follow the case of compact K\"ahler manifolds treated in
 \cite{DS8}. The fact that the functions in $\Fc$ are equal to a fixed smooth function near the boundary of $\B$ allows to adapt,  without difficulty, the integration by parts used in the compact setting, i.e. to estimate the mass of $d(u-u_n)\wedge d^c(u-u_n)\wedge T$. Note that we can extend the definition of $T\wedge S$ to the case where the potentials of $S$ are continuous outside a finite set of points. In this case, we need however some extra arguments.

We have the following result which was obtained in \cite{DS8}  when $X$ is a compact  K\"ahler manifold.

\begin{theorem} \label{th_wedge_ddc}
The $(p+1,p+1)$-current $T\wedge S$ is well-defined and is positive $\ddc$-closed. Moreover, it depends continuously on $T$ and on $S$ in the following sense. Let $T_n$ be positive $\ddc$-closed $(p,p)$-currents supported by a fixed closed set $F\subset X$ and converging to $T$. Let  $S_n$ be positive closed $(1,1)$-currents such that near $F$ we can write locally $S_n=\ddc u_n$ and $S=\ddc u$  with $u_n, u$ continuous p.s.h. and $u_n$ converging uniformly to $u$. Then $T_n\wedge S_n\to T\wedge S$. 
\end{theorem}

\section{Currents on compact K\"ahler manifolds} \label{section_Kahler}

In this section, we recall some classical results from the Hodge theory on compact K\"ahler manifolds and further properties of positive closed currents on such manifolds. We refer the readers to \cite{Demailly1,DS5, Voisin} for details.

\medskip
\noindent
{\bf $\bullet$ Hodge cohomology on compact K\"ahler manifolds.}
Let $(X,\omega)$ be a compact Hermitian manifold of dimension $k$. We say that the Hermitian form $\omega$ is {\it a K\"ahler form} if it is closed, i.e. $d\omega=0$. From now on, we assume that $\omega$ is a fixed K\"ahler metric on $X$. At each point $a\in X$, we can  find  local coordinates $z$ such that
$z=0$ at $a$ and $\omega$ is equal near 0 to 
$\ddc\|z\|^2$ modulo a term of order $\|z\|^2$. So, at the
infinitesimal level, a K{\"a}hler metric is close to the Euclidean
one. This is a crucial property in Hodge
theory on compact K\"ahler manifolds. Note that complex submanifolds of $X$ are also compact K\"ahler manifolds since $\omega$ restricted to these submanifolds defines a K\"ahler metric.

Recall that {\it the de Rham cohomology group} $H^l(X,\C)$ is
the quotient of the space of smooth closed $l$-forms by the
subspace of exact $l$-forms. This complex vector space is of finite dimension.
The real groups  $H^l(X,\R)$ are defined in the
same way using real forms and can be identified with a subspace of $H^l(X,\C)$. We also have 
$$H^l(X,\C)= H^l(X,\R)\otimes_\R\C.$$

In the definition of de Rham cohomology, we can also use currents instead of forms.
If $\alpha$ is a closed $l$-form or current, its class in $H^l(X,\C)$ is denoted
by $\{\alpha\}$. The group $H^0(X,\C)$ is just the set of constant
functions. So it is canonically identified to $\C$. The group $H^{2k}(X,\C)$ is also
isomorphic to $\C$. The isomorphism is given by the canonical map
$\{\alpha\}\mapsto\int_X\alpha$. 

For $l,m$ such that $l+m\leq 2k$, {\it the cup-product}
$$\smallsmile\ :\ H^l(X,\C)\times H^m(X,\C)\rightarrow H^{l+m}(X,\C)$$ 
is defined by 
$\{\alpha\}\smallsmile \{\beta\}:=\{\alpha\wedge\beta\}$. The Poincar\'e duality
theorem says that the cup-product is a non-degenerate bilinear form
when $l+m=2k$. So it defines an isomorphism between $H^l(X,\C)$ and
the dual of $H^{2k-l}(X,\C)$. 

Let $H^{p,q}(X,\C)$, $0\leq p,q\leq k$, denote the subspace of
$H^{p+q}(X,\C)$ generated by the classes
of closed $(p,q)$-forms or currents.
We call $H^{p,q}(X,\C)$ the {\it Hodge cohomology group}. 
Hodge theory shows that
$$H^{l}(X,\C)=\bigoplus_{p+q=l} H^{p,q}(X,\C) \qquad \mbox{and}\qquad
H^{q,p}(X,\C)=\overline{H^{p,q}(X,\C)}.$$ 
This, together with the Poincar{\'e} duality, induces a canonical isomorphism
between $H^{p,q}(X,\C)$ and the dual space of
$H^{k-p,k-q}(X,\C)$. 
Define for $p=q$
$$H^{p,p}(X,\R):=H^{p,p}(X,\C)\cap H^{2p}(X,\R).$$
We have
$$H^{p,p}(X,\C)=H^{p,p}(X,\R)\otimes_\R \C.$$
If $V$ is a complex subvariety of codimension $p$ in $X$, denote for simplicity $\{V\}$ the class in $H^{p,p}(X,\R)$ of the current $[V]$ of integration on $V$. 

Recall that {\it the Dolbeault cohomology group} $H^{p,q}_{\dbar}(X)$ is
the quotient of the space of $\dbar$-closed $(p,q)$-forms by the
subspace of $\dbar$-exact $(p,q)$-forms. Observe that a 
$(p,q)$-form is $d$-closed if and only if it is $\partial$-closed and
$\dbar$-closed.  By Hodge theory, we have the following natural isomorphism
$$H^{p,q}(X,\C)\simeq H^{p,q}_\dbar(X).$$
The result is a consequence of the following theorem, the so-called
{\it $\ddc$-lemma}, see e.g. \cite{Demailly3, Voisin}.

\begin{theorem} Let $\varphi$ be a smooth $d$-closed $(p,q)$-form or current on
  $X$. Then $\varphi$ is $\ddc$-exact if and only if it is $d$-exact (or 
  $\partial$-exact or $\dbar$-exact).
\end{theorem}

So $H^{p,q}(X,\C)$ is equal to the quotient of the space of $d$-closed $(p,q)$-forms by the subspace of $\ddc$-exact $(p,q)$-forms. If $T$ is a $\ddc$-closed $(p,p)$-current, it induces a linear form on $H^{k-p,k-p}(X,\C)$ and by Poincar\'e duality, it defines a class in $H^{p,p}(X,\C)$.  

The following result was obtained by Nguyen and the first author in \cite{DN1}. It generalizes previous results by Khovanskii, Teissier and Gromov \cite{Gromov1,Khovanskii,Teissier,Timorin}. 
\begin{theorem} \label{th_Dinh_Nguyen}
Let $\omega_1,\ldots,\omega_{k-1}$ be K\"ahler forms on $X$. Then the quadratic form $Q$ on $H^{1,1}(X,\R)$ given by
$$Q(c,c'):=-c\smallsmile c'\smallsmile \{\omega_1\}\smallsmile \cdots\smallsmile \{\omega_{k-2}\}$$
is positive definite on the hyperplane 
$$P:=\big\{c\in H^{1,1}(X,\R),\quad c\smallsmile \{\omega_1\}\smallsmile \cdots\smallsmile \{\omega_{k-1}\}=0\big\}.$$
\end{theorem}

When the forms $\omega_i$ are equal, we obtain the classical Hodge-Riemann theorem. Applying the above theorem to a class in the intersection of $P$ with the plane generated by $c,c'$, we get the following useful corollary.

\begin{corollary} \label{cor_HR}
Let $c$ and $c'$ be two linearly independent classes in $H^{1,1}(X,\R)$ such that $c^2\geq 0$ and $c'^2\geq 0$. Then $c\smallsmile c'\not=0$. In particular, if $c$ and $c'$ are two linearly independent classes in the closure $\overline \Kc$ of the K\"ahler cone, then $c\smallsmile c'\not=0$.
\end{corollary}

\medskip\noindent
{\bf $\bullet$ Projective manifolds.} An important large class of K\"ahler manifolds is the family of projective manifolds. They are isomorphic to complex submanifolds of projective spaces. We recall now this notion and fix some notation. 

The complex projective space $\P^k$ is a compact complex manifold of
dimension $k$. It
is obtained as the quotient of $\C^{k+1}\setminus\{0\}$ by the natural
multiplicative action of $\C^*$. In other words, $\P^k$ is the
parameter space of the complex lines through $0$ in
$\C^{k+1}$.  The image of a linear subspace
of dimension $p+1$ of $\C^{k+1}$ is a complex submanifold of dimension $p$ in
$\P^k$, biholomorphic to $\P^p$, and is
called {\it a projective subspace of dimension $p$}. {\it Hyperplanes}
of $\P^k$ are projective
subspaces of dimension $k-1$. The group $\GL(\C,k+1)$ of invertible
linear endomorphisms of $\C^{k+1}$ induces the group $\PGL(\C,k+1)$ of
automorphisms of $\P^k$ which is the projectivization of  $\GL(\C,k+1)$. It acts transitively on $\P^k$ and sends
projective subspaces to projective subspaces. 

Let $w=(w_0,\ldots,w_k)$ denote the standard coordinates
of $\C^{k+1}$. Consider the equivalence relation: {\it $w\sim w'$ if there
is $\lambda\in\C^*$ such that $w=\lambda w'$}. The projective space $\P^k$ is the
quotient of $\C^{k+1}\setminus\{0\}$ by this relation. We can
cover $\P^k$ by open sets $U_i$ associated with the open sets
$\{w_i\not=0\}$ in $\C^{k+1}\setminus \{0\}$. Each $U_i$ is
bi-holomorphic to $\C^k$ and $(w_0/w_i,\ldots,
w_{i-1}/w_i,w_{i+1}/w_i,\ldots,w_k/w_i)$ is a coordinate system on
this chart. The complement of $U_i$ is the hyperplane defined by $\{w_i=0\}$.
So $\P^k$ can be considered as 
a natural compactification of $\C^k$. We denote by $[w_0:\cdots:w_k]$
the point of $\P^k$ associated with $(w_0,\ldots, w_k)$. This expression
is {\it the homogeneous coordinates} on $\P^k$. 

The projective space $\P^k$ admits a K{\"a}hler form $\omega_\FS$,
called {\it the Fubini-Study form}. It 
is defined on the chart $U_i$ by 
$$\omega_\FS:=\ddc\log \Big(\sum_{j=0}^k \Big|{w_j\over
  w_i}\Big|^2\Big)^{1/2}.$$
In other words, if $\pi:\C^{k+1}\setminus\{0\}\rightarrow\P^k$ is the
canonical projection, then $\omega_\FS$ is defined by
$$\pi^*(\omega_\FS):=\ddc\log \Big(\sum_{i=0}^k \big|w_i|^2\Big)^{1/2}.$$

One can check that $\omega_\FS^k$ is a probability measure on
$\P^k$. 
The cohomology groups of $\P^k$ are very
simple. We have $H^{p,q}(\P^k,\C)=0$ for $p\not=q$ and
$H^{p,p}(\P^k,\C)\simeq\C$. The groups $H^{p,p}(\P^k,\R)$
and $H^{p,p}(\P^k,\C)$ are generated by the class of $\omega_\FS^p$.

Complex submanifolds of $\P^k$ are 
K\"ahler, as submanifolds of a K\"ahler manifold.
One has just to restrict the original K\"ahler form.
Chow's theorem says that such a manifold is 
{\it algebraic}, i.e. it is the set of common zeros of a finite family of
homogeneous polynomials in $w$. If $V$ is a subvariety of dimension $k-p$ of $\P^k$, we call {\it degree} of $V$ the number $\deg(V)$ of points in the intersection of $V$ with a generic projective subspace of dimension $p$. 
The class of the current of integration on $V$ is equal to $\deg(V)$ times the class of $\omega_\FS^p$.

\medskip\noindent
{\bf $\bullet$ Quasi-p.s.h. and d.s.h. functions.}
By maximum principle, p.s.h. functions on a compact connected manifold are
constant. However, a main interest of p.s.h. functions is their type of
singularities.  S.T. Yau introduced in \cite{Yau} the useful
notion of quasi-p.s.h. functions. 

\begin{definition} \rm
A {\it quasi-p.s.h. function} is
locally the difference of a p.s.h. function and a smooth one. A subset of $X$ is called {\it pluripolar} if it is contained in the pole set $\{u=-\infty\}$ of a quasi-p.s.h. function $u$. 
\end{definition}

Several properties of quasi-p.s.h. functions can be deduced from
properties of p.s.h. functions. We have the following result. Recall that $X$ is assumed to be compact.

\begin{proposition}
If $u$ is quasi-p.s.h. on $X$, it belongs to $L^p(X)$ for every $1\leq p<\infty$ and $\ddc u\geq -c\omega$ for some constant $c\geq 0$. 
If $(u_n)$ is a decreasing sequence of quasi-p.s.h. functions on $X$ satisfying $\ddc u_n\geq -c\omega$ with $c$ independent of $n$, then its limit is also a quasi-p.s.h. function.  If $S$ is a positive closed $(1,1)$-current and $\alpha$ a smooth real $(1,1)$-form in 
the cohomology class of $S$, then there is a quasi-p.s.h. function $u$, unique up to an additive constant, such that $\ddc u=S-\alpha$.  
\end{proposition}

The following regularization result is due to Demailly \cite{Demailly1}.

\begin{theorem}
Let $u$ be a quasi-p.s.h. function on $X$. Then there is a decreasing sequence of smooth functions $u_n$ with $\ddc u_n\geq -c\omega$ for some constant $c>0$ which converges pointwise to $u$.  In particular, if $S$ is a positive closed $(1,1)$-current on $X$ and $c>0$ is a constant large enough depending on $S$,  then the current $S+c\omega$ is limit of  smooth positive closed $(1,1)$-forms. 
\end{theorem}

\begin{example} \rm \label{ex_qpsh}
Let $S$ be a positive closed $(1,1)$-current of mass 1 on $\P^k$. Then there is a quasi-p.s.h. function $u$ on $\P^k$ such that $\ddc u=S-\omega_\FS$. In particular,  $S$ is cohomologous to $\omega_\FS$ and cannot be supported by a compact subset of the affine chart $\C^k\subset\P^k$ given by $\{w_0\not=0\}$. The restriction of $\omega_\FS$ to $\C^k$ is equal to $\ddc \log(1+\|z\|^2)^{1/2}$. Therefore, the function $u':=u+\log(1+\|z\|^2)^{1/2}$ is p.s.h. and is a potential of $S$ on $\C^k$, i.e. we have $S=\ddc u'$ on $\C^k$. Observe that the function $u'-\log^+\|z\|$ is bounded above. So we say that $u'$ is a {\it logarithmic potential} of $S$ on $\C^k$. We will use this example in Section \ref{section_Green}.
\end{example}

Recall also the notion of d.s.h. functions (differences of quasi-p.s.h. functions) which were introduced by
the authors  as observables in dynamics, see e.g. \cite{DS6}. They satisfy strong compactness properties and  are invariant under the action of holomorphic maps. 

\begin{definition} \rm   
A function on $X$ is called {\it d.s.h.} if it is equal outside a
pluripolar set to the difference of two quasi-p.s.h. functions. We
identify two d.s.h. functions if they are equal outside a pluripolar
set. 
\end{definition}

Let $\DSH(X)$ denote the space of d.s.h. functions on $X$. We
deduce easily from properties of p.s.h. functions the following result.

\begin{proposition}
The space $\DSH(X)$ is
contained in $L^p(X)$ for $1\leq p<\infty$. If $u$ is d.s.h. then we can write $\ddc
u=S^+-S^-$ where $S^\pm$ are positive closed $(1,1)$-currents in a same cohomology class. 
Conversely, if $S^\pm$ are positive
closed $(1,1)$-currents in a same cohomology class, then there is a
d.s.h. function $u$, unique up to a constant, such that $\ddc
u=S^+-S^-$. 
\end{proposition}

\begin{example} \rm
With the notation as in Example \ref{ex_qpsh}, the function $u'$ is d.s.h. on $\P^k$. We have $\ddc u'=S-[H]$ where $H:=\{w_0=0\}$ denotes the hyperplane at infinity.
\end{example}

We introduce a norm on $\DSH(\P^k)$. Let $S$ be a $(1,1)$-current which is the difference of two positive closed $(1,1)$-currents. Define
$$\|S\|_*:=\inf \|S^+\|+\|S^-\|,$$
where the infimum is taken over all positive closed $(1,1)$-currents $S^\pm$ such that $S=S^+-S^-$. 

Define for $u\in\DSH(X)$
$$\|u\|_\DSH:=|\langle \omega^k,u\rangle|+\|\ddc u\|_*.$$
The first term in the last definition
can be replaced with $\|u\|_{L^p}$,
 $1\leq p<\infty$; we then obtain equivalent norms. The space of d.s.h. functions endowed with
the above norm is a Banach space. We can also replace $\omega^k$ with any positive measure for which quasi-p.s.h. functions are integrable and get an equivalent norm. For example, we will need the following result.

\begin{lemma} \label{lemma_dsh_bound}
Let $\Fc$ be a family of d.s.h. functions on $X$ such that $\|\ddc u\|_*$ is bounded by a fixed constant for every  $u\in\Fc$. Let $U$ be a non-empty  open subset of $X$. Then $\Fc$ is bounded in $\DSH(X)$ if and only  if $\Fc$ is bounded in $L^1(U)$.
\end{lemma} 

The following result is also deduced from properties of p.s.h. functions.

\begin{theorem} \label{th_dsh_exp}
The canonical embedding of $\DSH(X)$ in $L^p(X)$ is compact. Moreover, if
$\Fc$ is a bounded subset of $\DSH(X)$, then there are positive constants $\alpha$ and $c$ such that 
$$\int_X e^{\alpha|u|} \omega^k\leq c.$$
\end{theorem}

We deduce the following useful corollary.

\begin{corollary} \label{cor_dsh_exp}
Let $\Fc$ be a bounded subset in $\DSH(X)$. If  $\nu$ is a probability measure associated with a bounded form of maximal degree on $X$, then there is a positive constant $c>0$ such that 
$$\langle \nu,|u|\rangle \leq c(1+\log^+\|\nu\|_\infty), \text{ where } \log^+:=\max(\log,0)$$
for all $u\in\Fc$. 
\end{corollary}
\proof
Let $\alpha$ be as in Theorem \ref{th_dsh_exp}. Using the concavity of the logarithm, we have 
$$\alpha \langle \nu,|u|\rangle = \big\langle \nu, \log e^{\alpha|u|} \big\rangle \leq \log  \big\langle \nu, e^{\alpha|u|} \big\rangle \leq \log\|\nu\|_\infty +\log \int_X e^{\alpha |u|}\omega^k.$$
Theorem \ref{th_dsh_exp} implies the result.
\endproof

\smallskip\noindent
{\bf $\bullet$ Positive $\ddc$-closed currents.} 
Consider now some properties of positive $\ddc$-closed $(p,p)$-currents. We have seen that they have cohomology classes in $H^{p,p}(X,\R)$. The following result was obtained in \cite{FS2}. 

\begin{theorem} \label{th_ddc_decom}
Let $S$ be a positive $\ddc$-closed $(1,1)$-current in $(X,\omega)$. Let $\alpha$ be a smooth real closed $(1,1)$-form in the class $\{S\}$. Then there is a $(0,1)$-current $\sigma$ such that 
$$S=\alpha+\partial \sigma+\overline{\partial \sigma}.$$
Moreover, the currents $\partial \overline \sigma, \dbar \sigma$ do not depend on the choice of $\sigma$ and are given by forms of class $L^2$. 
They vanish if and only if $S$ is closed.
\end{theorem}

We also need the following result.

\begin{theorem} \label{th_ddc_pot}
Let $S,\alpha$ and  $\sigma$ be as above. Assume that $S$ is smooth on an open set $U$. Then $\partial \overline \sigma$ and $\dbar \sigma$ are also smooth on $U$.  
Let $R$ be a positive closed $(k-2,k-2)$-current smooth outside a compact subset $K$ of $U$. Then 
 there is a constant $c>0$ depending only on $(X,\omega)$ such that 
$$\Big|\int_X \dbar \sigma\wedge \partial\overline \sigma\wedge R\Big|\leq c\|S\|^2\|R\|.$$
\end{theorem}
\proof
Observe that $\dbar \sigma\wedge \partial\overline \sigma$ is weakly positive.
The current $\partial\overline \sigma$ is the unique solution of the equation $\dbar(\partial\overline \sigma)=\partial S$ which is $\partial$-exact. To see this point, we can suppose $S=0$ and hence $\partial \overline\sigma$ is a holomorphic $(2,0)$-form whose cohomology class vanishes. It is well-known in Hodge theory that such a form on a compact K\"ahler manifold vanishes identically. Indeed, the group $H^{2,0}(X,\C)$ is isomorphic to the space of holomorphic $(2,0)$-forms.

Locally, the above equation admits a unique solution up to a holomorphic $(2,0)$-form which is smooth. So locally using classical explicit integral formula, we conclude that $\partial\overline \sigma$ is smooth on $U$. Its conjugate $\dbar \sigma$ is also smooth. 

The last estimate was obtained in \cite{FS2}. For smooths forms, it is obtained by integration by parts. So one has to use the regularization method in \cite{DS8} to find smooth positive $(k-2,k-2)$-forms $R_n$ on $X$ converging to a current $R'\geq R$ such that $\|R_n\|\leq c\|R\|$ for some constant $c\geq 0$ depending only on $(X,\omega)$. Moreover, $R'$ is smooth outside $K$ and $R_n$ converge to $R'$ locally uniformly on $X\setminus K$. 

Recall that $\dbar\sigma\wedge \partial\overline \sigma$ is a weakly positive $(2,2)$-form. So
using the smooth case applied to $R_n$, we obtain for a suitable constant $c$
$$\Big|\int_X \dbar \sigma\wedge \partial\overline \sigma\wedge R\Big|\leq \lim_{n\to\infty} \Big|\int_X \dbar \sigma\wedge \partial\overline \sigma\wedge R_n\Big|\leq c\|S\|^2\|R\|.$$
This completes the proof of the theorem.
\endproof

We also have the following version of Theorem \ref{th_ext_current}.

\begin{theorem} \label{th_ext_current_bis}
Let $E$ be an analytic subset of dimension smaller or equal to $k-p-1$ of a compact K\"ahler manifold $X$ of dimension $k$. 
Let $T$ be a positive $\ddc$-closed $(p,p)$-current on  $X\setminus E$. 
Then $T$ has finite mass and its extension by $0$ to $X$ is a positive $\ddc$-closed current.
\end{theorem}

When $X$ is an open manifold, the mass of $T$ is still locally bounded in $X$ but its extension by 0, still denoted by $T$, satisfies $\ddc T\leq 0$. On a compact K\"ahler manifold, Stokes' theorem implies that $\int_X\ddc T\wedge\omega^{k-p-1}=0$; it follows that $\ddc T=0$.

\medskip
\noindent
{\bf $\bullet$ Intersection number.} Let $T$ be a positive $(k-1,k-1)$-current on $X$ and $S$ a positive closed $(1,1)$-current on $X$. Assume that $\ddc T$ is a current of order 0 and that  the local potentials of $S$ are continuous on a neighbourhood of $\supp(T)$. Write $S=\alpha+\ddc u$ where $\alpha$ is a smooth $(1,1)$-form and $u$ a quasi-p.s.h. function continuous on a neighbourhood of $\supp(T)$. We define the intersection number $\langle T,S\rangle$ by
$$\langle T,S\rangle := \langle T,\alpha\rangle +\langle \ddc T,u\rangle.$$

It is not difficult to see that the definition does not depend on the choice of $\alpha$ and of $u$. Moreover, this number depends continuously on $T$ and on $S$ in the following sense. Let $T_n$ be a sequence of positive $(k-1,k-1)$-currents with support in a fixed compact set $K$ in $X$ such that $T_n\to T$ and $\|\ddc T_n\|$ are bounded uniformly on $n$. Let $S_n=\alpha_n+\ddc u_n$ be positive closed $(1,1)$-currents where $\alpha_n$ are smooth $(1,1)$-forms converging uniformly to $\alpha$ and $u_n$ quasi-p.s.h. functions continuous on a fixed neighbourhood $U$ of $K$ which converge to $u$ uniformly  on $U$. Then we have $\langle T_n,S_n\rangle \to \langle T,S\rangle$.

\section{Rigidity and Ahlfors-Nevanlinna currents} \label{section_rigid}

In this section, we introduce a notion of rigid sets and rigid cohomology classes. We also discuss some relations with the classical Ahlfors-Nevanlinna theory. We will see later that  Julia sets of H\'enon type maps are rigid and the cohomology classes of Green currents of automorphisms on compact K\"ahler surfaces satisfy a similar property. For the Ahlfors-Nevanlinna theory, we refer the readers to \cite{AS, Ne}. 

\begin{definition} \rm
Let $K$ be a closed subset of a complex manifold $X$. We say that $K$ is {\it $p$-rigid} (resp. {\it very $p$-rigid}) in $X$ if $K$ supports at most a non-zero positive closed (resp.  $\ddc$-closed) current of bidimension $(p,p)$ up to a multiplicative constant. The support $J$ of this current is called {\it the essential part} of $K$; by convention, $J$ is empty when such a current does not exist. For simplicity, when $p=1$, we say
 that $K$ is {\it rigid} (resp. {\it very rigid}) in $X$.
\end{definition}

Observe that if $K$ is $p$-rigid (resp. very $p$-rigid) as above and if $S$ is a non-zero positive closed (resp. $\ddc$-closed) current supported in $K$ then $S$ is {\it extremal}. That is, if $S=S_1+S_2$ with $S_1,S_2$ positive closed (resp. $\ddc$-closed) then $S_1,S_2$ are proportional to $S$. 
We also deduce that $J=\supp(S)$ is connected. 
Indeed, the restriction of $S$ to each connected component is also positive closed (resp. $\ddc$-closed).
Assume now that $(X,\omega)$ is a Hermitian manifold of dimension $k$. The following property is a direct consequence of the above definition.

\begin{proposition} \label{prop_rigid}
Let $K$ be a $p$-rigid (resp. very $p$-rigid) compact subset of $X$. Let $(\tau_n)$ be a sequence of positive currents of bidimension $(p,p)$ whose supports converge  to $K$. Assume that the masses $\|\tau_n\|$ converge to $1$ and  $d\tau_n\to 0$ (resp. $\ddc \tau_n\to 0$). Then $(\tau_n)$ converges to $\tau$ which is  the unique positive closed ($\ddc$-closed) current of mass $1$ and of bidegree $(p,p)$ supported on $K$.  
\end{proposition}
\proof
Since $K$ is compact and the masses $\|\tau_n\|$ are uniformly bounded, the sequence $(\tau_n)$ is relatively compact. If $\tau$ is a the limit of a subsequence, then it is a positive current of mass 1 supported on $K$ and by hypothesis, it is closed (resp. $\ddc$-closed). Since $K$ is $p$-rigid (resp. very $p$-rigid), $\tau$ is  the unique positive closed ($\ddc$-closed) current of mass $1$ supported on $K$.  
\endproof

In what follows, we only consider the case where $p=1$. Some results can be extended to higher dimensional currents but  in order to simplify the exposition we will discuss the topic only briefly in the last section.

\begin{definition} \rm
Let $\phi_n:\overline D_n\to X$ be a sequence of smooth maps on closed discs $\overline D_n$ in $\C$ which are holomorphic on $D_n$. Denote respectively by $a_n$ and $l_n$  the area of $\phi_n(D_n)$ and the length of $\phi_n(bD_n)$   counted with multiplicity. If $l_n=o(a_n)$ as $n\to\infty$, we say that $\phi_n:\overline D_n\to X$ is {\it an Ahlfors sequence of holomorphic discs}.
\end{definition}

With such an Ahlfors sequence, we can associate a sequence of currents of integration
$$\tau_n:={1\over a_n} (\phi_n)_*[D_n].$$
More precisely, if $\varphi$ is a smooth $(1,1)$-form on $X$ we have
$$\langle \tau_n,\varphi\rangle := {1\over a_n} \int_{D_n} \phi_n^*(\varphi).$$
Since $\|\tau_n\|=1$, see also Theorem \ref{th_wirtinger}, this sequence is relatively compact. 

\begin{definition} \rm
We call {\it Ahlfors current} any cluster value of the sequence $\tau_n$, i.e. the limit of a subsequence of $(\tau_n)$.
\end{definition}

This class of currents and the class of Nevanlinna's currents that will be defined later, were considered in several works, see e.g. \cite{BS,FS2,MQ}.
The condition $l_n=o(a_n)$ implies that $\|d\tau_n\|\to 0$. It follows that Ahlfors currents are positive closed and supported on the  cluster set of the sequence $\phi_n(D_n)$, i.e. on the set $\cap_{N\geq 0} \overline{\cup_{n\geq N} \phi_n(D_n)}$. 
We deduce from Proposition \ref{prop_rigid} that if all $\phi_n$ have images in a rigid compact set $K$ then $\tau_n$ 
converge to the unique positive closed current $\tau$ of mass 1 on $K$. Hence,
the essential part $J=\supp (\tau)$ of $K$ is contained in the cluster set of the sequence $\phi_n(\overline D_n)$. 

\medskip

We have the following classical result in Ahlfors-Nevanlinna theory.

\begin{theorem} \label{th_Nevanlinna}
Let $\phi:\C\to X$ be a non-constant holomorphic map into a Hermitian manifold $(X,\omega)$ with image in a compact set $K$. Let $\phi_r$ denote the restriction of $\phi$ to the disc $\overline \D_r$ of center $0$ and of radius $r$. Then there is a finite length subset $E\subset \R_+$ such that if $(r_n)\subset \R_+\setminus E$ and $r_n\to\infty$, the sequence $\phi_{r_n}$ is an Ahlfors sequence. In particular, if $K$ is rigid, its essential part $J$ is contained in $\overline{\phi(\C)}$. If $\phi$ has image in $J$ then $\phi(\C)$ is dense in $J$.
\end{theorem}
\proof
The first assertion is an interpretation of \cite{A,FS2}. We deduce that $K$ supports some Ahlfors currents of mass 1 with support in $\overline{\phi(\C)}$. So if $K$ is rigid, these Ahlfors currents are equal to the unique positive closed current of mass 1 on $K$. We deduce that $J$ is contained in  $\overline{\phi(\C)}$ and is equal to this set if $\phi$ has images in $J$. 
\endproof

Recall also the following useful re-parametrization lemma \cite{Ko}. 

\begin{lemma}[Brody-Zalcman] \label{lemma_Brody}
Let $h_n:\D\to X$ be a sequence of holomorphic maps from the unit disc to a complex manifold $X$ with images in a compact set $K$. Assume that the sequence $(h_n)$ is not locally equicontinuous.
Then there are an increasing sequence $(n_i)$ of positive integers, a sequence $(r_i)$ of positive numbers with $r_i\to\infty$ and a sequence of
affine maps $A_i:\D_{r_i}\to \D$ such that $h_{n_i}\circ A_i$ converge locally uniformly on $\C$ to a non-constant holomorphic map $\phi:\C\to X$ with image in $K$.
\end{lemma}

We deduce from the above results the following corollary.

\begin{corollary} \label{cor_Brody_rigid} 
With the notation as in the last lemma, if $K$ is rigid, then its essential part $J$ is contained in the cluster set of the sequence $h_n(\D)$, i.e. in the set $\cap_{N\geq 1} \overline{\cup_{n\geq N} h_n(\D)}$. 
\end{corollary}

We recall the notion of Kobayashi hyperbolicity on a complex manifold $X$ \cite{Ko}. Let $x$ be a point in $X$ and $\xi$ a complex tangent vector of $X$ at $a$. Consider the holomorphic maps $\phi:\D\to X$ on the unit disc in $\C$ such that $\phi(0)=x$ and $D\phi(0)=\lambda\xi$, where $D\phi$ denotes the differential of $\phi$ and $\lambda$ is a constant. 
{\it The Kobayashi-Royden infinitesimal pseudo-metric} is defined by 
$$\Kob_X(x,\xi):=\inf_\phi|\lambda|^{-1}.$$
It measures the size of discs that can be holomorphically  sent in $X$: the bigger is the disc, the smaller is the infinitesimal Kobayashi metric. If $X$ contains  non-constant holomorphic images of $\C$ passing through $x$ in the direction $\xi$, then $\Kob_X(x,\xi)=0$.

Kobayashi-Royden pseudo-metric is contracting for holomorphic maps: if $\psi:Y\to X$ is a holomorphic map between complex manifolds, then 
$$\Kob_X(\psi(y),D\psi(y)\zeta)\leq \Kob_Y(y,\zeta) \quad \text{for }  \zeta \text{ tangent to } Y \text{ at } y.$$
The Kobayashi pseudo-distance between two points is obtained by integrating the infinitesimal metric along curves between these points and taking the infimum. One obtain a metric on $X$ if  we have locally  $\Kob_X(x, \xi)\geq  c \|\xi\|$ for some constant $c>0$.

The Kobayashi-Royden infinitesimal pseudo-metric on $\D$ coincides with the Poincar\'e metric. A complex manifold $X$ is {\it Kobayashi hyperbolic} if $\Kob_X$ is a metric. In which case, holomorphic maps from any complex manifold $Y$ to $X$ are locally equicontinuous with respect to the Kobayashi metric on $X$. 

If $Z$ is a complex submanifold in $X$ or an open subset of  $X$, we say that $Z$ is {\it hyperbolically embedded} in $X$ if there is a continuous function $c>0$ on $X$ such that $\Kob_Z(z,\zeta)\geq c(z)\|\zeta\|$ for every complex tangent vector $\zeta$ of $Z$ at $z$. In particular, $Z$ is Kobayashi hyperbolic. Note that if $X$ is compact we can take $c$ constant.

\begin{proposition} \label{prop_non_hyp_emb}
Let $K$ be a compact subset of $X$ and $U$ a connected component of the interior of $K$. Assume that $U$ is not hyperbolically embedded in $X$. Then $\overline U$ contains a non-constant holomorphic image of $\C$. In particular, $\overline U$ supports an Ahlfors current. If $K$ is rigid, its essential part $J$ is contained in $\overline U$. 
\end{proposition}
\proof
Since $U$ is not hyperbolically embedded in $X$, there is a sequence of holomorphic maps $\phi_n:\D\to U$ such that $\|D\phi_n(0)\|\geq n$. Hence, the sequence $(\phi_n)$ is not locally equicontinuous on $\D$. By Lemma \ref{lemma_Brody}, there is a non-constant holomorphic map $\phi:\C\to \overline U$. The result then follows from Theorem \ref{th_Nevanlinna}. 
\endproof

Consider now a more general situation in Nevanlinna theory. Let $\Sigma$ be a Riemann surface. Assume that there is a function $\sigma:\Sigma \to [0,c_0[$ with $c_0\in[1,\infty]$ such that $\log\sigma$ is subharmonic and is harmonic out of a compact subset of $\Sigma$. Assume also that $\sigma$ is exhaustive, i.e. the set $\{\sigma\leq c\}$ is relatively compact in $\Sigma$ for every $c<c_0$. The main example in this work is the function $\sigma(\xi)=|\xi|$ for $\Sigma$ the unit disc $\D$ or the complex plane $\C$.  Note however that if $\Sigma$ is parabolic, i.e. an open Riemann surface without  non-constant bounded subharmonic functions, then $\Sigma$ admits a function $\sigma$ as above with $c_0=\infty$, see \cite{AS}. For example, the complement of a closed polar subset of $\C$ is parabolic.  

Let $\phi:\Sigma\to X$ be a non-constant holomorphic map with image in a compact set $K$ of a Hermitian manifold $(X,\omega)$. The main question in Nevanlinna theory is to study the value distribution of $\phi$. Define
$$\Sigma_t:=\big\{\xi\in\Sigma,\ \sigma(\xi)<t\big\} \quad \text{and} \quad A(t):=\int_{\Sigma_t} \phi^*(\omega).$$
So $A(t)$ is the area of $\phi(\Sigma_t)$ counted with multiplicity. The growth of $\phi$ is measured by the following Ahlfors-Nevanlinna function
$$T(r):=\int_0^r {A(t)\over t} dt=\int_\Sigma \log^+{r\over\sigma} \phi^*(\omega) \quad \text{with} \quad \log^+:=\max(\log,0).$$
So $T(r)$ is a weighted area of $\Sigma_r$. 

In what follows, we always make the hypothesis that $T(r)\to\infty$ when $r\to c_0$. It is not difficult to see that this property always holds when $\phi$ is non-constant and $c_0=\infty$, e.g. when $\Sigma=\C$ and $\sigma(\xi)=|\xi|$ or when $\Sigma$ is parabolic \cite{BS}. When $c_0$ is finite, the property is stronger than the fact that $A(t)\to \infty$ as $t\to c_0$. So it is not always true for $\Sigma=\D$ and $\sigma(\xi)=|\xi|$. However, it can be checked in some natural setting, e.g. when $\phi$ is the universal covering map of a generic leaf of a foliation by Riemann surfaces in $\P^2$, see \cite{FS2}. 

\begin{definition} \rm \label{def_Ne_current}
We call {\it Nevanlinna currents} any cluster value, for $r\to c_0$, of the family of currents
$$\tau_r:={1\over T(r)} \phi_*\Big(\log^+{r\over\sigma}\Big).$$
More precisely, if $\psi$ is a smooth $(1,1)$-form on $X$, we have
$$\langle \tau_r,\psi\rangle ={1\over T(r)} \int_\Sigma \Big(\log^+{r\over  \sigma}\Big) \phi^*(\psi)={1\over T(r)} \int_0^r{dt\over t} \int_{\Sigma_t} \phi^*(\psi).$$
\end{definition}

Nevanlinna currents describe the asymptotic behavior of the map $\phi$. The following result can be easily extended to a sequence of maps from Riemann surfaces to $X$ under appropriate conditions on the Nevanlinna characteristic functions. For simplicity, we only consider the case of one map $\phi:\Sigma\to X$.  

\begin{proposition}
With the above notation, assume also that $T(r)\to \infty$ when $r\to c_0$. Then all Nevanlinna currents associated with $\phi$ are positive $\ddc$-closed currents of mass $1$ with support in $\overline{\phi(\Sigma)}$. 
\end{proposition}
\proof
By definition, the currents $\tau_r$ are positive, of mass 1 and supported on $\overline{\phi(\Sigma)}$ which is contained in the compact set $K$. The last assertion in the proposition is clear. So we only have to check that Nevanlinna currents are $\ddc$-closed. We have 
$$\ddc \tau_r={1\over T(r)} \phi_*\Big(\ddc \log^+{r\over\sigma}\Big).$$

Since $\log\sigma$ is harmonic out of a compact set, for $r$ close enough to $c_0$, the expression in the last parenthesis is equal to the difference of two positive measures with compact support $\nu_r-\nu$, where $\nu_r$ is supported on the real curve $\{\sigma=r\}$ and $\nu:=\ddc\log\sigma$. Since $\log^+{r\over\sigma}$ vanishes out of $\Sigma_r$, by Stokes' formula we have 
$$\|\nu_r\|-\|\nu\|=\langle \nu_r-\nu,1\rangle =\Big\langle \ddc \log^+{r\over\sigma}, 1\Big\rangle=
\Big\langle \log^+{r\over\sigma}, \ddc 1\Big\rangle=0.$$
So the mass of $\nu_r$ is independent of $r$. It follows that $\|\ddc\tau_r\|\to 0$ since $T(r)\to\infty$.   
We deduce that Nevanlinna currents are $\ddc$-closed. This completes the proof of the proposition.

Note that when $c_0$ is finite or when $\Sigma=\C$ and $\phi(\C)$ not contained in a compact curve one can show that all Nevanlinna currents are supported in  $\cap_{r<c_0}\overline{\phi(\Sigma\setminus \Sigma_r)}$.
\endproof

Note that in comparison with Theorem \ref{th_Nevanlinna} there is here no exceptional set $E$: any cluster value of $(\tau_r)$ is $\ddc$-closed. 
As above, we obtain the following result.

\begin{corollary}
With the above notation, assume that $T(r)\to\infty$ and  that the compact set $K$ is very rigid. Then all Nevanlinna currents on $K$ are equal to the unique positive $\ddc$-closed current of mass $1$ on $K$ and the essential part $J$ of $K$ is contained in  $\overline{\phi(\Sigma)}$. If $\phi$ has images in $J$ then $\phi(\Sigma)$ is dense in $J$.
\end{corollary}

We now introduce a notion of rigidity for cohomology classes which also appears naturally  in dynamics on compact K\"ahler manifolds. 

\begin{definition} \rm
Let $X$ be a compact K\"ahler manifold. Let $c$ be a class in $H^{p,p}(X,\R)$. We say that $c$ is {\it rigid} (resp. {\it very rigid}) if it contains a non-zero positive closed $(p,p)$-current $T$ which is the unique positive closed (resp. $\ddc$-closed) $(p,p)$-current in the class $c$. 
\end{definition}

So if $(T_n)$ is a sequence of positive closed (resp. $\ddc$-closed) $(p,p)$-currents such that the sequence of classes $\{T_n\}$ converges to a rigid (resp. very rigid) class $c$ as above, then $T_n$ converge to $T$ as currents.

\section{Basic properties of H\'enon type maps on $\C^2$} \label{section_Henon}

In this section, we give some elementary properties of the Fatou-Julia theory in the setting of H\'enon type maps in $\C^2$. We first describe the group of polynomial automorphisms of $\C^2$. The main references to this section are \cite{BS1,FS4,FM,Sibony}.

Let $p$ be a polynomial in one complex variable of degree $d\geq 2$. A {\it H\'enon map} in $\C^2$ has the following form 
$$h(z_1,z_2):=\big(p(z_1)+az_2,z_1\big),$$
where $a$ is a constant in $\C^*$. It is clear that $h$ is an automorphism of $\C^2$, i.e. a bijective holomorphic map on $\C^2$; the inverse map is given by
$$h^{-1}(z_1,z_2)=\Big(z_2,{z_1-p(z_2)\over a}\Big).$$
We will see later that 
these automorphisms have rich dynamics. 

There are also the so-called {\it elementary automorphisms}. They have the following form
$$e(z_1,z_2):=\big(az_1+p(z_2),bz_2+c\big)$$
where $a,b,c$ are constants in $\C$ with $a,b\not=0$,  and $p$ a polynomial of degree $d\geq 0$. The automorphism  $e$ preserves the fibration $\{z_2=\const\}$. Their dynamics is simple to analyze. The family of elementary automorphisms, that we denote by $E$,  has a group structure. We have the following important result \cite{Jung}.

\begin{theorem}[Jung]
The group $\Aut(\C^2)$  of polynomial automorphisms of $\C^2$ is generated by the elementary ones and the group $A$ of complex affine automorphisms. More precisely, $\Aut(\C^2)$ is the amalgamated product of $A$ and $E$ along their intersection $A\cap E$. 
\end{theorem}

Friedland and Milnor deduced from Jung's theorem the following property \cite{FM}.

\begin{corollary}[Friedland-Milnor] \label{cor_FM}
Let $f$ be a polynomial automorphism of $\C^2$. Then $f$ is conjugate in $\Aut(\C^2)$ either to an elementary automorphism or to a H\'enon type map, i.e. a finite composition of H\'enon maps. 
\end{corollary}

So for the dynamical study of a polynomial automorphism,  we can restrict our attention to polynomial automorphisms of the second type, i.e. finite compositions of $(p_j(z_1)+a_jz_2,z_1)$ with $p_j$ polynomials of degree $d_j\geq 2$ and 
$a_j\in\C^*$. 
It is quite fruitful to consider the extension of such an automorphism as a birational map in the projective plane $\P^2$. In this paper we favor  this point of view. Let us be more precise. In $\P^2$ we consider the homogeneous coordinates $[w_0:w_1:w_2]$. The line at infinity is defined by the equation $w_0=0$. 

Given any polynomial map $g:\C^2\to\C^2$ of the form $g=(P_1,P_2)$ with $\max(\deg P_1,\deg P_2)=d\geq 1$, we can consider the extension $\overline g$ to $\P^2$:
$$\overline g[w]:=\Big[w_0^d:w_0^dP_1\big({w_1\over w_0},{w_2\over w_0}\big):w_0^dP_2\big({w_1\over w_0},{w_2\over w_0}\big)\Big].$$
If we restrict to $\C^2$, i.e. to the affine chart $\{w_0\not=0\}$ with affine coordinates $(z_1,z_2):=(w_1/w_0,w_2/w_0)$, we find the original map. Note that $\overline g$ is not defined on a finite set at infinity where all three coordinate functions of $\overline g$ vanish. We call this set {\it the indeterminacy set} of $\overline g$. It is non-empty in general.

For example, if $h$ is a H\'enon map as above, then
$$\overline h[w_0:w_1:w_2]=\Big[w_0^d:w_0^dp\big({w_1\over w_0}\big)+aw_2w_0^{d-1}:w_1w_0^{d-1}\Big].$$
We observe that the three coordinate functions vanish simultaneously at only a point $I_+=[0:0:1]$. So $I_+$ is the unique indeterminacy point of $\overline h$. Similarly the indeterminacy point for $\overline h^{-1}$ is $I_-=[0:1:0]$. Observe that $\{w_0=0\}\setminus I_+$ is mapped by $\overline h$ to $I_-$ and $\{w_0=0\}\setminus I_-$ is contracted by $\overline h^{-1}$ to $I_+$. 

It is useful to note that $I_-$ is fixed for $\overline h$ and that it is super-attracting i.e. both eigenvalues of the differential of $\overline h$ vanish at this point. Indeed if we express $\overline h$ in the chart $\{w_1\not= 0\}$ with affine coordinates $z_0':=w_0/w_1, z_2':=w_2/w_1$, we get 
$$\overline h(z_0',z_2')=\Big({{z_0'}^d\over {z_0'}^dp({1\over z_0'})+az'_2{z'_0}^{d-1}}, {{z_0'}^{d-1}\over {z_0'}^dp({1\over z_0'})+az_2'{z_0'}^{d-1}}\Big)$$ 
and it is then easy to compute the eigenvalues of the Jacobian matrix at $I_-=(0,0)$. The point $I_+$ is also fixed and super-attracting for $\overline h^{-1}$. 
We deduce from this discussion the following properties. 

\begin{proposition} \label{prop_Henon_type}
Let $f$ be a H\'enon type map on $\C^2$. For simplicity, denote also by $f$ its extension as a birational map on $\P^2$. Then the indeterminacy set of $f$ (resp. $f^{-1}$) is reduced to the point $I_+=[0:0:1]$ (resp. $I_-=[0:1:0]$). Moreover, $f$ (resp. $f^{-1}$) contracts $\{w_0=0\}\setminus I_+$  (resp. $\{w_0=0\}\setminus I_-$) to the fixed super-attracting point $I_-$ (resp. $I_+$). 
\end{proposition}

Because of Corollary \ref{cor_FM}, we will be interested only in the dynamics of H\'enon type maps $f$  as above. 
Denote by $\Gamma$ the closure of the graph of $f:\P^2\setminus I_+\to\P^2$ in $\P^2\times\P^2$. Since $f$ is a birational map, $\Gamma$ is an irreducible algebraic set of dimension 2. 
The image of $\Gamma$ by the involution $(z,z')\mapsto (z',z)$ on $\P^2\times\P^2$ is the closure of the graph of $f^{-1}$.
Let $\pi_1,\pi_2$ denote the canonical projections from $\P^2\times\P^2$ onto its factors. The following lemma gives a description of $\Gamma$.

\begin{lemma} \label{lemma_Henon_graph}
We have 
$$\pi_1^{-1}(z)\cap \Gamma=(z,f(z)) \quad \text{for}\quad  z\in\P^2\setminus I_+,
\qquad
\pi_1^{-1}(I_+)\cap\Gamma= I_+\times \{w_0=0\}$$ 
and
$$\pi_2^{-1}(z)\cap \Gamma=(f^{-1}(z),z) \quad \text{for}\quad  z\in\P^2\setminus I_-, \qquad
\pi_2^{-1}(I_-)\cap\Gamma= \{w_0=0\}\times I_-.$$ 
\end{lemma}
\proof
The first identity follows from the definition of $\Gamma$. 
By continuity, $\Gamma$ contains the points $(f^{-1}(z),z)$ with $z\in\P^2\setminus I_-$ since this is true for $z\in\C^2$. 
So the third identity is obtained similarly. 
We prove the second one and the rest can be obtained in the same way. 
By Proposition \ref{prop_Henon_type}, $\{w_0=0\}\setminus I_-$ is sent by $f^{-1}$ to $I_+$. It follows that 
$\pi_1^{-1}(I_+)\cap\Gamma\supset I_+\times \{w_0=0\}$. The third identity implies that $\pi_1^{-1}(I_+)$ contains no point in $I_+\times \C^2$. The result follows.
\endproof

Note that $f^n$ is also a H\'enon type map for $n\geq 1$. So Proposition \ref{prop_Henon_type} and Lemma \ref{lemma_Henon_graph} apply to $f^n$. 

Define $\widetilde U_+$ as the basin of attraction of $I_-$ for $f$ in $\P^2\setminus I_+$, $U_+:=\widetilde U_+\cap \C^2$ and $K_+:=\C^2\setminus U_+$. So $\widetilde U_+,U_+$ are open sets in $\P^2$ and $K_+$ is closed in $\C^2$. The following lemma implies that the similar objects associated with $f^n$, $n\geq 1$, do not depend on $n$.

\begin{proposition} \label{prop_K_J}
We have $\widetilde U_+\setminus U_+=\{w_0=0\}\setminus I_+$.
The closure of $K_+$ in $\P^2$ is $\overline K_+=K_+\cup I_+$.  
We also have that
$$K_+=\big\{z\in\C^2,\ (f^n(z))_{n\geq 0} \text{ is bounded in } \C^2\big\}.$$ 
Moreover, $U_+,K_+$ and $\partial K_+$ are invariant under $f$ and under $f^{-1}$.
\end{proposition}
\proof
By definition, $f^{\pm 1}(U_+)\subset  U_+$. Since $f$ is an automorphism of $\C^2$, we deduce that $f(U_+)=f^{-1}(U_+)=U_+$. So $U_+$ is invariant under $f^{\pm 1}$. We deduce that $K_+$ and then $\partial K_+$ are also invariant under $f^{\pm 1}$.
 
By Proposition \ref{prop_Henon_type}, $\{w_0=0\}\setminus I_+$ is contained in $\widetilde U_+$. It follows that $\overline K_+\subset K_+\cup I_+$. If $\overline K_+\not= K_+\cup I_+$, there would be a ball $B$ centered at $I_+$ such that $\partial B\subset \widetilde U_+$. If $n$ is large enough, $f^n(\partial B)$ is close to $I_-$. The classical Hartogs theorem for domains in  $\C^2$ implies that $f^n$ extends to a holomorphic map on $B$ and contradicts the fact that $I_+$ is an indeterminacy point for $f^n$. So $\overline K_+=K_+\cup I_+$.

It remains to show that $K_+$ is the set of points of bounded orbit in $\C^2$. It is clear that such a point belongs to $K_+$. Let $z$ be a point in $K_+$. We have to show that the sequence $(f^n(z))_{n\geq 0}$ is bounded in $\C^2$. 

Since $I_+$ is fixed super-attracting for $f^{-1}$, there is a neighbourhood $W$ of $I_+$ such that $f^{-1}(W)\Subset W$. We can choose $W$ small enough so that $z\not\in W$. The property of $W$ implies that if $f^n(z)\not\in W$ then 
$f^{n+1}(z)\not\in W$.  Since $K_+$ is invariant, we deduce that $f^n(z)$ is in  $K_+\setminus W$ which is a compact set in $\C^2$. This completes the proof of the proposition.
\endproof

\begin{definition} \rm
We call {\it Fatou set} of $f$ the largest open set $F_+$ in $\P^2\setminus I_+$ on which the sequence $(f^n)_{n\geq 1}$ is locally equicontinuous as maps from $F_+$ to $\P^2$. The {\it Julia set} of $f$ is the complement $J_+$ of $F_+$ in $\P^2$. 
\end{definition}

The following result is a consequence of the last proposition.

\begin{corollary}
We have $J_+=\partial \overline K_+$. The Fatou set $F_+$ is the union of $\widetilde U_+$ and the interior of $K_+$. 
\end{corollary}
\proof
If $z$ is a point in $\partial K_+$, then any neighbourhood of $z$ contains both points with bounded orbits and points with orbits going to $I_-$. So $(f^n)_{n\geq 0}$ cannot be equicontinuous on any neighbourhood of $z$. We deduce that $\partial \overline K_+\subset J_+$. Since the complement of $\partial \overline K_+$ is the union of $\widetilde U_+$ and the interior of $K_+$, it is enough to check that this set is contained in $F_+$.

Since $\widetilde U_+$ is the basin of $I_-$, it is contained in $F_+$. Consider a small open set $U$ in $K_+$. We only have to show that $(f^n)_{n\geq 0}$ is locally equicontinuous on $U$. Choose $W$ as in the proof of Proposition \ref{prop_K_J} small enough so that $W\cap U=\varnothing$. So $f^n$ restricted to $U$ has images in the compact set $K_+\setminus W$ for every $n$. We deduce using Cauchy's formula for holomorphic functions that $(f^n)_{n\geq 0}$ is locally equicontinuous on $U$. The result follows.
\endproof

We can define the Julia and Fatou sets for $f^{-1}$ in the same way. Let $\widetilde U_-$ denote the basin of $I_+$ for the map $f^{-1}$ on $\P^2\setminus I_-$. Define $U_-:=\widetilde U_-\cap\C^2$ and $K_-:=\C^2\setminus U_-$. We obtain as above the following results.

\begin{proposition} \label{prop_K_J_bis}
We have $\widetilde U_-\setminus U_-=\{w_0=0\}\setminus I_-$. 
The closure of $K_-$ in $\P^2$ is $\overline K_-=K_-\cup I_-$. 
We also have
$$K_-=\big\{z\in\C^2,\ (f^{-n}(z))_{n\geq 0} \text{ is bounded in } \C^2\big\}.$$ 
 Moreover, $U_-,K_-$ and $\partial K_-$ are invariant under $f$ and under $f^{-1}$. The Julia set $J_-$ of $f^{-1}$ is equal to $\partial \overline K_-$ and the Fatou set $F_-$ is the union of $\widetilde U_-$ and the interior of $K_-$. 
\end{proposition}

Define 
$$K:=K_+\cap K_-.$$
This set is compact in $\C^2$ and is invariant under $f$ and under $f^{-1}$.
 
Denote by $a$ the determinant of the complex Jacobian matrix of $f$ with respect to the coordinates $(z_1,z_2)$. Since $f$ is a polynomial automorphism, $a$ is an invertible  polynomial. So it is constant.
We have the following result.

\begin{proposition} \label{prop_Henon_Jac}
If $|a|>1$ then $K_+$ has empty interior, i.e. $J_+=\partial \overline K_+=\overline K_+$. If $|a|<1$ then $K_-$ has empty interior, i.e. $J_-=\partial \overline K_-=\overline K_-$. If $|a|=1$, then $K_+\setminus K$ and $K_-\setminus K$ have empty interior and all connected components of the interior of $K$ are periodic.
\end{proposition}
\proof
We consider the case $|a|\geq 1$. The case $|a|\leq 1$ can be treated in the same way. Let $W$ be as in Proposition \ref{prop_K_J}. Consider the sequence of sets
$$H_n:=f^n(K_+\setminus (K\cup W))=K_+\setminus (K_-\cup f^n(W)).$$
Since $\widetilde U_-$ is the basin of $I_-$ for $f^{-1}$, the sequence $f^n(W)$ increases to $\widetilde U_-$. So the sequence $(H_n)$ decreases to the empty set.

On the other hand, the Euclidean volume of $H_n$, which is compact in $\C^2$,  is larger or equal to the Euclidean volume of $H_0$ because $|a|\geq 1$. Therefore, all these sets have zero volume. We deduce that $K_+\setminus K$ has empty interior and zero volume. If $|a|>1$, since $K$ is invariant, the same argument implies that $K$ has zero volume and empty interior. We conclude that $J_+=\partial \overline K_+=\overline K_+$.

Assume now that $|a|=1$. The above discussion can be applied to $f$ and $f^{-1}$. So  $K_+\setminus K$ and $K_-\setminus K$ have empty interior and zero volume. If $U$ is a connected component of the interior of $K$, then $f^n(U)$ is also a connected component of the interior of $K$ with the same volume. Since $K$ has finite volume, there are integers $n<m$ such that $f^n(U)=f^m(U)$. It follows that $f^{m-n}(U)=U$. So $U$ is periodic. This completes the proof of the proposition.
\endproof

We will need the following notion.

\begin{definition} \rm
Let $g=(P,Q)$ be a map from $\C^2$ to $\C^2$, where $P$ and $Q$ are polynomials on $\C^2$. We call {\it algebraic degree} of $g$ the maximum of the degrees of $P$ and of $Q$. We denote it by $\adeg(g)$ or $\deg(g)$ if there is no risk of confusion.
\end{definition}

\begin{proposition} \label{prop_deg}
Let $f$ be a H\'enon type map on $\C^2$ of algebraic degree $d$. Then the algebraic degree of $f^n$ is equal to $d^{|n|}$ for $n\in\Z$. Moreover, the first (resp. second) coordinate function of $f^n$ for $n\geq 0$ (resp. $n\leq 0$) is of degree $d^{|n|}$ and its homogeneous part of maximal degree is a monomial in $z_1$ (resp. $z_2$). The second (resp. first) coordinate function of $f^n$ for $n\geq 0$ (resp. $n\leq 0$) is of degree $d^{|n|-1}$.
\end{proposition}
\proof
Let $h_i(z)=(p_i(z_1)+a_iz_2,z_1)$ be H\'enon maps of algebraic degrees $d_i:=\deg p_i\geq 2$ with $a_i\in\C^*$ and $1\leq i\leq m$. We first show that the algebraic degree of $g:=h_1\circ\cdots\circ h_m$ is equal to $d_1\ldots d_m$. 
It is clear that this degree is smaller or equal to $d_1\ldots d_m$. We easily obtain by induction on $m$ that the degrees  of the coordinate functions of $g$ are $d_1\ldots d_m$ and $d_2\ldots d_m$ respectively and the homogeneous part of maximal degree of the first coordinate function of $g$ is a monomial in $z_1$. So the algebraic degree of $g$ is $d_1\ldots d_m$. 

Applying this property to $f^n$ with $n\geq 0$, we obtain the proposition for $n\geq 0$. The case $n\leq 0$ can be obtained by induction on $m$ as above. It is enough to use the form of $h^{-1}$ given above and to observe that $\deg(h)=\deg(h^{-1})$.
\endproof

A more conceptual way to look at the algebraic degree is as follows. Let $V$ be an algebraic curve in $\P^2$ which does not contain $I_-$. If $\Gamma$, $\pi_1$ and $\pi_2$ are as above, then $f^{-1}(V)=\pi_1(\pi_2^{-1}(V)\cap\Gamma)$. We see that $f^{-1}(V)$ depends continuously on $V$ under the hypothesis that $I_-\not\in V$. If $V$ is given in $\C^2$ by $\{P=0\}$ where $P$ is a polynomial in $z_1,z_2$ such that $\deg (P)=\deg(V)$, then $f^{-1}(V)$ is given by $\{P\circ f=0\}$. For generic $V$, we see that $\deg(f^{-1}(V))=d\deg(V)$ and by continuity the property holds for all $V$ with $I_-\not\in V$. We conclude that the action of $f^*$ on $H^{1,1}(\P^2,\R)$ is just the multiplication by $d$. 

Now if $f'$ is another H\'enon type map of degree $d'$ then since $f^{-1}(V)$ does not contain $I_-$, we have $\deg (f\circ f')^{-1}(V)=dd'\deg(V)$ and therefore $\deg(f\circ f')=\deg(f)\deg(f')$. This also allows us to obtain the first assertion of the last proposition. The idea here can be used in the more general setting of algebraically stable meromorphic maps considered  in \cite{FS1,Sibony}.

We close this section by observing that since periodic points in $\C^2$ belong to $K$, the following result shows that $K$ is always non-empty.

\begin{proposition}
Let $f$ be a H\'enon type automorphism of degree $d$ of $\C^2$. Then $f$ admits $d^n$ periodic points of period $n$ in $\C^2$ counted with multiplicity. Moreover, it admits an infinite number of distinct periodic orbits and has periodic points of arbitrary large period.
\end{proposition}
\proof
A theorem by Shub-Sullivan \cite[p.323]{KH} implies that if $a$ is an isolated fixed point for all $f^n$ then the multiplicity of $a$ as a fixed point of $f^n$ is bounded uniformly on $n$.  Therefore, the second assertion in the proposition is a consequence of the first one.

For the first assertion, since $f^n$ is a H\'enon type map of degree $d^n$, we only need to consider the case $n=1$. If $a$ is a fixed point in $\C^2$ we can associate it to the point $(a,a)$ in the intersection of $\Gamma$ with the diagonal $\Delta$ of $\P^2\times\P^2$. The above description of $f$ at infinity shows that outside $\C^2\times\C^2$ the graph $\Gamma$ intersects $\Delta$ transversally at two points $I_+\times I_+$ and at $I_-\times I_-$. So the number of fixed points in $\C^2$ counted with multiplicity is equal to $\{\Gamma\}\smallsmile \{\Delta\}-2$. 

It follows from K\"unneth formula \cite{Voisin}  that $\{\Delta\}$ is equal to $\{a\times\P^2\}+\{\P^2\times a\}+\{L\times L\}$ where $a$ is a point  and $L$ is a projective line in $\P^2$. Choose $a$ and $L$ generic. We see that  $\{a\times\P^2\}\smallsmile \{\Gamma\}=1$ since $a\times \P^2$ intersects $\Gamma$ transversally at a point. We also have $\{\P^2\times a\}\smallsmile \{\Gamma\}=1$. The intersection of $L\times L$ with $\Gamma$ is the set of points $(a,f(a))$ with $a\in L$ and $f(a)\in L$. So $\{L\times L\}\smallsmile \{\Gamma\}$ is the number of points in $L\cap f^{-1}(L)$. By B\'ezout theorem, this number is equal to $d$ since $\deg(L)=1$ and $\deg(f^{-1}(L))=d$.  The result follows.
\endproof

\section{Green currents and rigidity of Julia sets} \label{section_Green}

In this section, we give a construction of the Green currents associated with H\'enon type maps and prove a rigidity property of the supports of these currents. This rigidity property is the key point in the proof of several dynamical properties of H\'enon type maps. For previous works see \cite{BS1, FS4, Sibony} and the references therein.

Let $f$ be a H\'enon type map on $\C^2$ as above with algebraic degree $d\geq 2$. Define for $n\geq 0$
$$G^+_n(z):={1\over d^n} \log^+\|f^n(z)\| \quad \text{with } \log^+:=\max(\log, 0).$$
These functions measure the convergence speed of orbits to infinity. We have the following theorem.

\begin{theorem} \label{th_Green}
The sequence $(G^+_n)_{n\geq 0}$ converges locally uniformly on $\C^2$
to a positive H\"older continuous p.s.h. function $G^+$. The convergence is almost decreasing: for every neighbourhood $W$ of $I_+$, there is a sequence of real numbers $(c_n)$ decreasing to $0$ such that $G^+_n+c_n$ decrease to $G^+$ on $\C^2\setminus W$. We have $G^+\circ f=dG^+$ on $\C^2$ and $G^+=0$ exactly on $K_+$. Moreover, $G^+$ is pluriharmonic and strictly positive on $U_+$ and $G^+(z)-\log^+|z_1|$ extends to a pluriharmonic function on a neighbourhood of $\{w_0=0\}\setminus I_+$ in $\P^2$. 
\end{theorem}

We need the following special case of \cite[Prop. 2.4]{DS1}.

\begin{lemma} \label{lemma_Holder}
Let $X$ be a metric space with finite diameter and $\Lambda:X\to X$ a Lipschitz map with $\dist(\Lambda(a),\Lambda(b))\leq A\dist(a,b)$ and $A>1$. Let $v$ be a bounded Lipschitz function on $X$. Then the series $\sum_{n\geq 0} d^{-n} v\circ \Lambda^n$ converge pointwise to a function which is $\beta$-H\"older continuous for any $\beta$ such that $0\leq \beta\leq 1$ and $\beta<{\log d\over \log A}\cdot$  
\end{lemma}

\noindent
{\bf Proof of Theorem \ref{th_Green}.} Since $I_+$ is attracting for $f^{-1}$ replacing $W$ with a suitable smaller
open set allows to assume that $f(X)\subset X$ for $X:=\P^2\setminus \overline W$. Define 
$$v(z):=G_1^+(z)-G_0^+(z)={1\over d} \log^+\|f(z)\|-\log^+\|z\|.$$
This function is Lipschitz on $X$. 
Observe that for $z\in X$ with $\|z\|$ large enough, $|z_2|$ is bounded by a constant times $|z_1|$.
Proposition \ref{prop_deg} for $n=1$ shows that $|v|$ is bounded on $X$ by a constant $c>0$. 

We have 
$$|G_{n+1}^+-G_{n}^+|=d^{-n} |v\circ f^{n}|\leq cd^{-n} \quad \mbox{on}\quad  X.$$
Therefore, $G_n^+$ converge uniformly on $X$ to a function $G^+$. So $G^+_n$ converge locally uniformly on $\C^2$ to $G^+$. Since the functions $G_n^+$ are positive p.s.h., $G^+$ is also positive p.s.h. Since $G_n^+\circ f=dG_{n+1}^+$, we have $G^+\circ f=dG^+$. 

Observe also that 
$$G^+(z)=\log^+\|z\|+\sum_{n\geq 0} d^{-n} v(f^n(z)).$$
Lemma \ref{lemma_Holder} implies that $G^+$ is H\"older continuous out of any given neighbourhood $W$ of $I_+$.
Moreover, we can take for the H\"older exponent any constant $\beta$, $0\leq\beta\leq 1$, which is  smaller than 
$${\log d\over \max_{\P^2\setminus W} \log\|Df\|}\cdot$$
By definition, the Green function of $f^n$ is also equal to $G^+$. Therefore, we can take for the H\"older exponent any constant $\beta$, $0\leq\beta\leq 1$,  smaller than 
$$\beta_0=\sup_{n\geq 1}{\log d\over \max_{\P^2\setminus W} \log\|Df^n\|^{1/n}}=\lim_{n\to\infty}{\log d\over \max_{\P^2\setminus W} \log\|Df^n\|^{1/n}}\cdot$$
Finally, since any compact set in $\P^2\setminus I_+$ is sent by some $f^m$ into $\P^2\setminus W$, we conclude that $\beta_0$ does not depend on the choice of the small open set $W$. 
Moreover, if $V$ is a neighbourhood of $\overline K_-$, we have 
$$\beta_0=\lim_{n\to\infty}{\log d\over \max_{V\setminus W} \log\|Df^n\|^{1/n}}\cdot$$
So $G^+$ is locally $\beta$-H\"older continuous on $\C^2\setminus I_+$ for every $\beta$ such that $0\leq\beta\leq 1$ and $\beta<\beta_0$.

The function $v$ is bounded, so the above identity relating $G^+$ and $v$  shows that $G^+>0$ in a neighbourhood $U$ of $I_-$. The invariance relation $G^+\circ f=dG_+$ implies that $G^+>0$ on the open set $\cup_{n\geq 0} f^{-n}(U)\cap \C^2$ which is equal to the basin  $U_+$ of $I_-$ in $\C^2$.
By Proposition \ref{prop_K_J}, if $z$ is in $K_+$, the orbit of $z$ is bounded in $\C^2$. Therefore, by definition of $G_n^+$, we have $G^+(z)=0$. So $G^+=0$ exactly on $K_+$ which is the complement of $U_+$ in $\C^2$. 

Define $c_n:=c(d^{-n}+d^{-n-1}+\cdots)$. Clearly, the sequence  $(c_n)$ decreases to 0. We have
$$(G_{n+1}^++c_{n+1})-(G_n^++c_n)=d^{-n} v\circ f^n-cd^{-n}\leq 0.$$
It follows that $G^+_n+c_n$ decrease to $G^+$ on $X$. 

Consider now the open set $X_R:=X\cap\{|z_1|>R\}$ with a constant $R$ large enough. Observe that $|z_2|$ is bounded by a constant times $|z_1|$ on $X_R$ and Proposition \ref{prop_deg} for $n=1$ implies that $f(X_R)\subset X_R$ and hence $f^n(X_R)\subset X_R$ for every $n\geq 0$. Denote by $f^n_1$ and $f^n_2$ the coordinate functions of $f^n$ and 
$$H^+_n:=d^{-n}\log |f^n_2| \quad \text{on} \quad X_R\cap \C^2.$$
Since $f_1^n=O(f_2^n)$ on $X_R$, we deduce that $G^+_n-H^+_n$ converge uniformly to 0. It follows that $H^+_n$ converge uniformly to $G^+$ on $X_R\cap \C^2$. 

Finally, since $H^+_n$ is pluriharmonic, $G^+$ is pluriharmonic on $X_R\cap\C^2$. The invariance relation of $G^+$ implies that it is pluriharmonic on $\cup_{n\geq 0} f^{-n}(X_R)\cap \C^2$ which is equal to $U_+$. Moreover, the function $H^+_n(z)-\log|z_1|$ extends to a pluriharmonic function through the line at infinity. We deduce that $G^+(z)-\log|z_1|$ is pluriharmonic in a neighbourhood of $\{w_0=0\}\setminus I_+$ in $\P^2$. This completes the proof of the theorem.
\hfill $\square$

\medskip

Define for each constant $c\geq 0$
$$K_+^c:=\{G^+\leq c\},  \  J_+^c:=\partial \overline K_+^c (=\overline{\{G^+=c\}} \text{ if } c>0)
\  \text{and} \  G^+_c:=\max(G^+-c,0).$$
Note that $\{G^+=c\}$, for $c>0$, is a real analytic hypersurface in $\C^2$. 
We have the following lemma.

\begin{lemma} \label{lemma_K_c}
For $c\geq 0$, we have $\overline K_+^c=K_+^c\cup I_+$ and  $J_+^c=\{G^+=c\}\cup I_+$. Moreover, $T_+^c:=\ddc G_c^+$ is a positive closed $(1,1)$-current on $\C^2$ which extends by zero to a positive closed $(1,1)$-current of mass $1$ supported on $J_+^c$. We also have $\supp(T^c_+)=J^c_+$ and the Hausdorff dimension of $J^c_+$ is strictly larger than $2$ on any open set which intersects $J^c_+$. The set $J_+^c$ is a smooth analytic real hypersurface of $\C^2$ for $c>0$. 
\end{lemma}
\proof
It is clear that $G_c^+$ is p.s.h. So $T_+^c$ is a positive closed $(1,1)$-current. Since $G_c^+$ is pluriharmonic on $\{G^+<c\}\cup \{G^+>c\}$, this current $T_+^c$ vanishes outside $\{G^+=c\}$. The maximum principle applied to $-G^+_c$ shows that $G^+_c$ is not pluriharmonic on any open set which intersects $J_+^c$. So $T_+^c$ does not vanish identically and its support is equal to $J^c_+$. We also deduce that $J^c_+$ is not compact in $\C^2$, see Example \ref{ex_qpsh}. The assertion on the Hausdorff dimension of $J^c_+$ is a consequence of the H\"older continuity of  the potential $G_c^+$ of $T_+^c$. Indeed, $T_+^c$ has no mass on sets of Hausdorff dimension $2+\epsilon$ for some constant $\epsilon>0$ small enough, see \cite[1.7.3]{Sibony}. Indeed, $\epsilon$ can be chosen equal to the above constant $\beta_0$ involved in the H\"older continuity of $G^+$. 

The property of $G^+(z)-\log|z_1|$ given in Theorem \ref{th_Green} implies that $\overline K^c_+\subset K_+^c\cup I_+$. Since $K^c_+$ contains the support of $T_+^c$ which is unbounded in $\C^2$, we conclude that $\overline K^c_+= K_+^c\cup I_+$ and then
$J_+^c=\{G^+=c\}\cup I_+$.  

We have $\overline K_+^c=K_+^c\cup I_+$ and $T_+^c$ is a positive closed current on $\P^2\setminus I_+$. 
By Theorem \ref{th_ext_current}, $T_+^c$ extends by 0 to a positive closed $(1,1)$-current on $\P^2$ that we still denote by $T_+^c$. We check that its mass is equal to 1.

Let $\widetilde L$ be a projective line through $I_-$. The intersection $T_+^c\wedge [\widetilde L]$ is a well-defined positive measure with compact support in the complex line $L:=\widetilde L\cap \C^2$. We only have to show that this measure has mass 1.  In the complex line $L$, it is equal to the function $\ddc {G^+_c}_{|L}$. Since the function ${G^+_c}_{|L}$  is subharmonic, harmonic near infinity and has logarithmic growth, it is well-know that the associated measure is a probability measure, see Example \ref{ex_qpsh} with $k=1$. This implies that the mass of $T_+^c$ is equal to 1. 

It remains to verify that $J_+^c$ is smooth for $c>0$. Since the family of these hypersurfaces is invariant under the automorphism $f$, it is enough to prove the property in a neighbourhood of $I_-$ for $c$ large enough.
Using the local coordinates $z_1':=1/z_1$ and $z_2'=z_2/z_1$, we see that  $J_+^c$ are level sets of a function of the form $\exp(-G^+)=\alpha |z_1'|(1+o(1))$ with $\alpha\not=0$. Clearly, these hypersurfaces are smooth near $I_-$. This
completes the proof of the lemma.
\endproof

\begin{definition} \rm
We call $G^+$ {\it the Green function} and $T_+$ {\it the Green current} associated with $f$.
\end{definition}

Here is a main result in this section.

\begin{theorem} \label{th_K_rigid}
The sets  $\overline K_+$ and $J_+$ are very rigid and the Green current $T_+$ is the unique positive $\ddc$-closed current of mass $1$ supported on these sets.
\end{theorem}

We will give later the proof of this theorem. We first show the following related result.

\begin{theorem} \label{th_equi_Green}
Let $U$ be a neighbourhood of $I_-$ and $A>0$ a constant. Let $(S_n)_{n\geq 1}$ be a sequence of positive closed $(1,1)$-currents of mass $1$ on $\P^2$. Assume that $S_n$ admits a quasi-potential $u_n$ such that $|u_n|\leq A$ on $U$ for every $n$. Then $d^{-n} (f^n)^*(S_n)$ converge to $T_+$ exponentially fast: there is a constant $c>0$ such that 
$$|\langle d^{-n}(f^n)^*(S_n)-T_+,\varphi\rangle|\leq cnd^{-n}\|\varphi\|_{\Cc^2}$$
for every test $(1,1)$-form $ \varphi$ of class $\Cc^2$. 
\end{theorem}
\proof
The idea is to observe that $u_n$ have bounded DSH norms. We will use Corollary \ref{cor_dsh_exp} after estimating the uniform norm of the $(2,2)$-form $(f^{-n})^*(\ddc\varphi)$ away from $I_-$.  

We can replace $U$ by a suitable smaller domain in order to assume that $f(U)\Subset  U$. 
Multiplying $\varphi$ with a constant allows to assume that $\|\ddc\varphi\|_\infty\leq 1$. So $\nu:=\ddc\varphi$ is  a complex measure of mass $\leq 1$.
Define $\nu_n:=(f^n)_*(\nu)$. Since $\nu$ has no mass at infinity and $f^n$ is an automorphism on $\C^2$, the measures $\nu_n$ and $\nu$  have the same mass which is smaller or equal to 1. 

Let $\nu_n'$ and $\nu_n''$ be the restrictions of $\nu_n$ to $\P^2\setminus U$ and to $U$ respectively. We have $\nu_n=\nu_n'+\nu_n''$, $\|\nu_n'\|\leq 1$ and $\|\nu_n''\|\leq 1$. Observe that $f^{-1}$ defines a map from $\P^2\setminus U$ to $\P^2\setminus U$ with $\Cc^1$-norm bounded by some constant $M$. So the $\Cc^1$-norm of $f^{-n}$ on $\P^2\setminus U$ is bounded by $M^{n}$. Since $\nu_n'$ is the pull-back of $\nu$ by this map, we have $\|\nu_n'\|_\infty\leq M^{4n}$. 

Let $g^+:=G^+-\log(1+\|z\|^2)^{1/2}$ be a quasi-potential of $T_+$. Theorem \ref{th_Green} implies that $g^+$ is smooth near $I_-$. Define $v_n:=u_n-g^+$. Replacing $A$ by a suitable constant allows to assume that $|v_n|\leq A$ on $U$.  
We have since $T_+$ is invariant
\begin{eqnarray*}
\langle d^{-n}(f^n)^*(S_n)-T_+,\varphi\rangle & = & \langle d^{-n}(f^n)^*(S_n)-d^{-n}(f^n)^*(T_+),\varphi\rangle \\
& = & \langle d^{-n} (f^n)^*(\ddc v_n),\varphi\rangle =d^{-n}\langle v_n, \ddc (\varphi\circ f^{-n})\rangle\\
& = & d^{-n} \langle v_n,(f^n)_*(\ddc\varphi)\rangle = d^{-n} \langle v_n,\nu_n\rangle\\
&=& d^{-n} \langle \nu_n',v_n\rangle+d^{-n} \langle \nu_n'',v_n\rangle.
\end{eqnarray*}

The second term in the last sum is of order $O(d^{-n})$ since $\|\nu''_n\|\leq 1$ and $|v_n|\leq A$ on the support of $\nu_n''$. By Lemma \ref{lemma_dsh_bound}, $v_n$ has bounded DSH-norm. Since $\|\nu_n'\|\leq 1$ and $\|\nu_n'\|_\infty\leq M^{4n}$, Corollary \ref{cor_dsh_exp} implies that the first  term in the last sum satisfies
$$|d^{-n}\langle \nu'_n,v_n\rangle|\leq d^{-n}c(1+\log^+ M^{4n})=O(nd^{-n}).$$
This completes the proof of the theorem. Note that one can easily deduce from the obtained estimate a similar estimate for H\"older continuous test forms using the classical interpolation theory between Banach spaces.  
\endproof

The following result shows for generic analytic subsets $V$ in $\P^2$ that  the sequence $f^{-n}(V)$ is asymptotically  equidistributed with respect to the current $T_+$, see also  \cite{BS0,FS4}. 

\begin{corollary}
Let $S$ be a positive closed $(1,1)$-current of mass $1$ on $\P^2$ whose support does not contain $I_-$. Then $d^{-n} (f^n)^*(S)\to T_+$ exponentially fast as $n\to\infty$. In particular, if 
$V$ is an analytic set of pure dimension $1$ which does not contain $I_-$, then  $d^{-n}[f^{-n}(V)]\to \deg(V)T_+$ exponentially fast.
\end{corollary}
\proof
The first assertion is a direct consequence of Theorem \ref{th_equi_Green}. In order to obtain the second assertion, 
it is enough to take $S:={1\over \deg(V)} [V]$ and to observe that $(f^n)^*(S)={1\over \deg(V)}[f^{-n}(V)]$. 
\endproof

\noindent
{\bf Proof of Theorem \ref{th_K_rigid}. }
Let $S$ be a positive closed $(1,1)$-current of mass 1 with support in $\overline K_+$. We first show that $S=T_+$. Define $S_n:=d^n(f^n)_*(S)$ on $\C^2$. This is a positive closed $(1,1)$-current on $\C^2$ with support in $K_+$. So it is also a positive closed $(1,1)$-current on $\P^2\setminus I_+$. By Theorem \ref{th_ext_current}, it extends by 0 through $I_+$ to a positive closed $(1,1)$-current on $\P^2$.  

We have $(f^n)^*(f^n)_*(S)=S$ on $\C^2$. Hence $d^{-n}(f^n)^*(S_n)=S$ on $\C^2$. It follows that $d^{-n}(f^n)^*(S_n)=S$ on $\P^2$, in particular, we have $\|S_n\|=1$ because $(f^n)^*$ multiplies the mass of a positive closed $(1,1)$-current by $d^n$. Since the currents $S_n$ vanish on a neighbourhood of $I_-$, Theorem \ref{th_equi_Green} implies that $S=T_+$.

Consider now a positive $\ddc$-closed current $S$ of mass 1 on $\overline K_+$. We have to prove that $S=T_+$. It is enough to show that $S$ is closed. We only have to check that $\partial S=0$ since we can obtain in the same way that $\dbar S=0$. 

Define $S_n$ as above. This current is positive $\ddc$-closed on $\P^2\setminus I_+$. By Theorem \ref{th_ext_current_bis}, this current is positive $\ddc$-closed on $\P^2$. We also deduce as above that $S=d^{-n}(f^n)^*(S_n)$ and $\|S_n\|=1$. By Theorems \ref{th_ddc_decom} and \ref{th_ddc_pot}, we can write 
$$S_n=\omega_\FS+\partial\sigma_n+\overline{\partial\sigma_n},$$
where $\sigma_n$ is a $(0,1)$-current such that  $\dbar\sigma_n$ and $\partial\overline\sigma_n$ are $L^2$ forms with norms bounded independently of $n$. 

We have
$$\partial S=d^{-n}(f^n)^*(\partial S_n)=-d^{-n}\dbar (f^n)^*(\partial\overline\sigma_n).$$
On the other hand, we have 
$$(f^n)^*(\dbar\sigma_n)\wedge (f^n)^*(\partial \overline\sigma_n)=(f^n)^*(\dbar \sigma_n\wedge\partial \overline\sigma_n).$$
The maximal degree form $\dbar \sigma_n\wedge\partial \overline\sigma_n$ defines a positive measure of finite mass and this mass is invariant under $(f^n)^*$. Since $\partial\overline\sigma_n$ is of bidegree $(2,0)$, we deduce that $d^{-n}(f^n)^*(\partial\overline\sigma_n)\to 0$ in $L^2$. Hence $\partial S=0$. This completes the proof of the theorem. 
\hfill $\square$

\begin{corollary}
Let $(S_n)_{n\geq 1}$ be a sequence of positive $(1,1)$-currents with support in a given compact set in $\P^2\setminus I_-$. Assume that $\|\ddc S_n\|=o(d^n)$ and that $\langle T_-,S_n\rangle\to1$. Then $d^{-n}(f^n)^*(S_n)\to T_+$ in $\P^2$. In particular, if 
$S$ is a positive $(1,1)$-current with compact support in $\P^2\setminus I_-$ such that $\ddc S$ is of order $0$, then $d^{-n}(f^n)^*(S)\to c_ST_+$ in $\P^2$, where $c_S:=\langle T_-,S\rangle$, see the end of Section \ref{section_Kahler}  for the definition of this intersection number.
\end{corollary}
\proof
We first show that the mass of $R_n:=d^{-n}(f^n)^*(S_n)$ converges to $1$. We have
$$\|R_n\|:=\langle d^{-n}(f^n)^*(S_n),\omega_\FS\rangle = \langle S_n, d^{-n} (f^n)_*(\omega_\FS)\rangle.$$
Theorem \ref{th_Green} applied to $f^{-1}$ implies that $d^{-n} (f^n)_*(\omega_\FS)$ admit quasi-potentials $v_n$ which decrease to a quasi-potential $v$ of $T_-$. Moreover, all these functions are continuous out of $I_-$ and we deduce from the proof of Theorem \ref{th_Green} that $v_n-v$ is bounded by a constant times $d^{-n}$ on each compact subset of $\P^2\setminus I_-$. So  we can write using the intersection numbers introduced in Section \ref{section_Kahler} and the convergence result established there,
\begin{eqnarray*}
\langle S_n, d^{-n} (f^n)_*(\omega_\FS)\rangle & = & \langle S_n, \omega_\FS +\ddc v_n\rangle=\langle S_n, T_-\rangle  +\langle S_n, \ddc (v_n-v)\rangle \\
& = & \langle S_n, T_-\rangle  +\langle \ddc S_n, v_n-v \rangle.
\end{eqnarray*}
By hypothesis, $\|\ddc S_n\|=o(d^n)$ and $\langle S_n,T_-\rangle\to 1$, hence the last sum converges to 1. Therefore, we get  $\|R_n\|\to 1$. 

Let $T$ be any cluster value of $R_n=d^{-n}(f^n)^*(S_n)$. This is a positive current of mass $1$. We have to show that it is equal to $T_+$.
Since $S_n$ have supports in a given compact subset of  $\P^2\setminus I_-$, the current $T$ is supported on $\overline K_+$. By Theorem \ref{th_K_rigid}, it is enough to check that $T$ is $\ddc$-closed or equivalently $\ddc R_n\to 0$. 

We have
$$\ddc R_n=d^{-n} (f^n)^*(\ddc S_n).$$
By hypothesis, $\|\ddc S_n\|=o(d^n)$. Since the operator $(f^n)^*$ preserves the mass of measures, we easily deduce that $\ddc R_n\to 0$. This completes the proof of the corollary. 
\endproof

\begin{remark} \rm
In the last result, the condition on the support of $S_n$ can be weakened. We will not try to get the most general statement here.  
Let $S$ be a positive $(1,1)$-current on $\P^2$ such that $\ddc S$ is of order $0$ with compact support in $\P^2\setminus I_-$. Define $c_S:=\langle T_-,S\rangle$ and denote by $c_S'$ the mass of $T_-\wedge S$ at $I_-$. 
We can show that  $d^{-n}(f^n)^*(S)\to (c_S-c_S')T_+$ in $\C^2$. 
The property applies to positive $\ddc$-closed currents, see also \cite{CG} for the case of positive closed currents.
\end{remark}

The following result gives us a property of the Fatou set as a consequence of the above rigidity properties. It can be applied to the basin of an attractive fixed point which is biholomorphic to $\C^2$. 

\begin{corollary} \label{cor_Fatou_Henon}
Let $\Omega$ be a Fatou component contained in $K_+$, i.e. a connected component of the interior of $K_+$. Then one of the following properties holds 
\begin{enumerate} 
\item $\Omega$ is not hyperbolically embedded in $\P^2$ and  $\partial \Omega= J_+$;
\item $\Omega$ is a connected component of the interior of $K$ and is periodic; in particular, it is bounded and hyperbolically embedded in $\C^2$;
\item $\Omega$ is hyperbolically embedded in $\P^2$; it  is a connected component of the interior of $K_+\setminus K_-$ and is periodic; in particular, $f^{-n}$ converges locally uniformly on $\Omega$ to $I_+$ as $n\to\infty$; 
\item $\Omega$ is hyperbolically embedded in $\P^2$ and is wandering, i.e. the open sets $f^n(\Omega)$, $n\in\Z$, are disjoint.
\end{enumerate}
\end{corollary}
\proof
If $\Omega$ is not hyperbolically embedded in $\P^2$, Proposition \ref{prop_non_hyp_emb} implies that $\overline\Omega$ contains $J_+$. 
Since $\Omega$ is a component of the interior of $K_+$, we deduce that $\partial\Omega=J_+$. 
Assume now that $\Omega$ is hyperbolically embedded in $\P^2$, non-wandering and is not contained in $K$. So $\Omega$ is periodic. By Proposition \ref{prop_K_J_bis}, we only have to check that $\Omega\cap K_-=\varnothing$. 

Assume that $\Omega\cap K_-\not =\varnothing$. Replacing $f$ by an iterate, we can assume that $\Omega$ is invariant.
Let $\phi:\D\to \Omega$ be a holomorphic disc centered at a point in $\Omega\cap K_-$ but not contained in $K_-$. 
So $\phi(\D)$ contains both points of bounded negative orbit and points of unbounded negative orbit. Hence the family of maps $\phi_n:=f^{-n}\circ \phi:\D\to\Omega$ is not locally equicontinuous on $\D$.  It follows that $\Omega$ is not hyperbolically embedded in $\P^2$. This is a contradiction.
\endproof

We also have the following result.

\begin{proposition} \label{th_Kc_rigid}
Assume that $K_+$ is not contained in a (possibly singular) real analytic hypersurface of a neighbourhood of $K_+$. Then for every $c\geq 0$, the set $\partial \overline K_+^c$ is very rigid and $T_+^c$ is the unique positive $\ddc$-closed $(1,1)$-current of mass $1$ supported on this set. 
\end{proposition}
\proof We only have to consider the case $c>0$. Recall that by Lemma \ref{lemma_K_c}, we have $\overline K_+^c=K_+^c\cup I_+$. Let $S$ be a positive $\ddc$-closed current of mass 1 on $\partial \overline K_+^c$. We have to show that $S=T_+^c$. 
We can define as in Theorem \ref{th_K_rigid} a positive $\ddc$-closed current $S_n$ of mass 1 on $\partial \overline K_+^{cd^n}$ such that $S=d^{-n}(f^n)^*(S_n)$. We obtain  as in this theorem that $S$ is closed. 

Let $u_n$ be a d.s.h. function such that $\ddc u_n=S_n-T_+^{cd^n}$. This function is pluriharmonic out of $\partial \overline K_+^{cd^n}$. Subtracting from $u_n$ a constant allows to assume  that
$u_n(a)=0$ at a given fixed point $a$. We easily deduce from the pluriharmonicity on $\{G^+<c\}$  that
the sequence $(u_n)$ is relatively compact in the space of d.s.h. functions. In particular, these functions are uniformly bounded in any compact subset of $\{G^+<c\}$. 

On the other hand, we have $\ddc (d^{-n} u_n\circ f^n)=S-T_+^c$. Therefore, $d^{-n} u_n\circ f^n$ is equal to $u_0$ plus a constant. The condition $u_n(a)=0$ implies that $d^{-n} u_n\circ f^n=u_0$. For any point $z\in K_+$, since $(f^n(z))_{n\geq 0}$ is relatively compact in $\{G^+<c\}$, we deduce from the last identity that $u_0(z)=0$. So $u_0=0$ on $K_+$ and hence $u_n=0$ on $K_+$. 
We show that $u_0=0$ on $\{G^+<c\}$. Assume this is not the case. Then $u_n$ does not vanish identically on $\{G^+<cd^n\}$. 
In particular, $K_+$ is contained in a real analytic hypersurface of $\{G^+<cd^n\}$, a contradiction.

Let $\widetilde L$ be a projective line through $I_-$. The function $u_0+G_c^+$, which is a potential of $S$ in $\C^2$,  is subharmonic on the complex line $L:=\widetilde L\setminus I_-$ with logarithmic growth. It vanishes on the open set $L\cap \{G^+<c\}$ and harmonic outside the real analytic curve $\{G^+=c\}\cap L$. It is clear that such a function should be the Green function associated with $\{G^+\leq c\}\cap L$. The function $G_c^+$ satisfies the same property on $L$. We conclude that $u_0+G^+=G^+$ and that $u_0=0$ on $L$. This property holds for all $L$. So we have $u_0=0$ on $\P^2$. The proposition follows.
\endproof

\begin{remarks} \rm
It is likely that the hypothesis on $K_+$ is always satisfied, see \cite{FS4} for results in this direction.
The above results can be applied to $f^{-1}, K_-,J_-, G^-$ and $T_-$. The invariant probability measure $\mu:=T_+\wedge T_-$ turns out to be the unique measure of maximal entropy $\log d$. It is exponentially mixing and saddle periodic points are equidistributed with respect to $\mu$. We refer to \cite{BLS2,Dinh1} for details.
\end{remarks}

Consider a family of H\'enon type maps $(f_c)_{c\in\Sigma}$ of degree $d$ depending holomorphically on a parameter $c$ in a complex manifold $\Sigma$. The maps $f_c$ have the same indeterminacy points $I_+=[0:0:1]$ and $I_-=[0:1:0]$. Denote by $K_\pm(c), J_\pm(c),G^\pm(c), T_\pm(c)$ and $\mu(c)$ the dynamical objects associated with $f_c$ constructed as above. 

We associate  this family with the dynamical system
$$F:\Sigma\times\P^2\to \Sigma\times \P^2 \quad \text{with} \quad F(c,z):=(c,f_c(z)).$$
The associated Green functions are defined by 
$$\Gc^\pm(c,z):=\lim_{n\to\infty} d^{-n} \log \|F^{\pm n}(c,z)\|.$$ 
We have seen that the above limit always exists. Define also
$$\Kc_\pm:=\big\{(c,z)\in\Sigma\times \C^2,\ (F^{\pm n}(c,z)) \text{ bounded in } \Sigma\times\C^2\big\}.$$ 
This is the union of the sets $\{c\}\times K_\pm(c)$ with $c\in\Sigma$. 
One can identify $K_\pm(c)$, $G^\pm(c),T_\pm(c)$ and $\mu(c)$ with the restriction to $\{c\}\times \C^2$ of $\Kc_\pm, \Gc^\pm, \ddc \Gc^\pm$ and $\ddc\Gc^+\wedge \ddc \Gc^-$ respectively. 

A point $(c,z)$ is {\it stable} if the sequence $(F^n)_{n\geq 0}$ is equicontinuous in a neigbourhood of $(c,z)$. 
As in the case without parameters, we see that the set of unstable points is exactly the boundary of $\Kc_+$ which is also the support of $\Tc_+:=\ddc \Gc^+$. 
We have the following theorem.

\begin{theorem}
The functions $\Gc^\pm$ are  locally H\"older continuous, p.s.h. positive and vanishes exactly on the set $\Kc_\pm$. The sets $\Kc_\pm$ are rigid and $\Tc_\pm:=\ddc  \Gc^\pm$ are the unique positive closed $(1,1)$-currents, up to a multiplicative constant, supported on $\Kc_\pm$. The supports of these currents are the boundaries of $\Kc_\pm$. Moreover, we have 
$$\partial \Kc_\pm=\overline{\bigcup_{c\in\Sigma} \{c\}\times \partial K_\pm(c)}.$$ 
\end{theorem}
\proof
The first assertion and the assertion on the supports of $\Tc_\pm$ are obtained as in the case without parameters. We have seen that $K_\pm(c)$ are rigid. Since $T_\pm(c)$ are slices of $\Tc_\pm$ with respect to the natural projection on $\Sigma$, we easily deduce from Theorem \ref{th_slice} that $\Kc_\pm$ are rigid. We prove now the last assertion for $\Kc_+$ and the case of $\Kc_-$ can be obtained in the same way.

First, since the support of $T_+(c)$, which is a slice of $\Tc_+$, is $\{c\}\times \partial K_+(c)$, we have the inclusion 
$$\overline{\bigcup_{c\in\Sigma} \{c\}\times \partial K_+(c)}\subset \partial \Kc_+.$$ 
Let $(c_0,z_0)$ be a point in $\partial \Kc_+$. Since this point is on the support of $\Tc_+$, the function $\Gc^+$ on $\Uc_+$ does not have a pluriharmonic extension on any neighbourhood of $(c_0,z_0)$. Fix a bidisc $D$ in $\C^2$ centered at $z_0$. By Levi's extension theorem, already used in Theorem \ref{th_slice}, there is $c$ arbitrary close to $c_0$ such that the restriction of $\Gc^+$ to $\{c\}\times D$ is not pluriharmonic. It follows that $\{c\}\times D$ intersects $\{c\}\times \partial K_+(c)$. Therefore, $(c_0,z_0)$ belongs to the set on the left hand side of the above inclusion. This completes the proof of the theorem.
\endproof

\begin{remarks} \rm
If we consider the maps $c\mapsto J_\pm(c)$, they are lower semi-continuous with respect to the Hausdorff metric for sets. We can look for the upper semi-continuous envelops, i.e. the maps whose graphs are the smallest closed sets containing the graphs of  $c\mapsto J_\pm(c)$. According to the previous theorem, these closed sets are equal to the boundaries of $\Kc_\pm$.  So the map $c\mapsto J_+(c)$ is continuous at a point $c_0$ if and only if the restriction of $\partial \Kc_+$ to $\{c_0\}\times  \C^2$ is equal to $\{c_0\} \times J_+(c_0)$. In particular, if the interior of $K_+(c_0)$ is made only of basins of attractive periodic cycles, then $c_0$ is a point of continuity. Indeed, it is not difficult to see that these basins are contained in the interior of $\Kc_+$. 

We know that $\Tc_+\wedge \Tc_-$ is a positive closed $(2,2)$-current. Its slices can be identified to the equilibrium measures $\mu(c)$ of $f_c$. The map $c\mapsto \supp\mu(c)$ is lower semi-continuous. The support of $\Tc_+\wedge \Tc_-$ is the graph of an upper semi-continuous map. Il would be interesting to show that the later map is the envelop of $c\mapsto \supp\mu(c)$ and to study the continuity points of these maps. 
\end{remarks}

\begin{remark} \rm
We can show that the current $\Tc_+$ is laminar. It is enough to apply the criteria obtained in \cite{dD,Dinh1}.
\end{remark}

\section{Automorphisms of  compact K\"ahler surfaces} \label{section_surface}

In this section we discuss automorphisms of positive entropy on a compact K\"ahler surface. We will see that several techniques presented in the last section can be adapted here. The fact that there are no  indeterminacy points simplifies the analytic part. In contrast, the actions of automorphisms on cohomology is simple but non-trivial.

Let $f:X\to X$ be a holomorphic automorphism on a compact K\"ahler surface $(X,\omega)$. 
The pull-back and push-forward actions of $f$ on differential forms induce linear operators $f^*$ and $f_*$ on $H^{p,q}(X,\C)$ with $0\leq p,q\leq 2$. Since $f$ is an automorphism, we have 
$$f^*=(f^{-1})_*, \quad f_*=(f^{-1})^* \quad \mbox{and} \quad f^*\circ f_*=f_*\circ f^*=\id.$$ 
Moreover,  we have $f^*=f_*=\id$ when $p=q=0$ or $p=q=2$.

\begin{definition} \rm
The spectral radius $d$ of $f^*$ on $H^{1,1}(X,\C)$ is called {\it the dynamical degree} of $f$. 
\end{definition}

It is easy to show that 
$$d=\lim_{n\to\infty} \Big(\int_X (f^n)^*\omega\wedge\omega\Big)^{1/n}= \lim_{n\to\infty} \Big(\int_X \omega\wedge (f^n)_*\omega\Big)^{1/n}.$$
We then deduce that the dynamical degree of $f^{-1}$ is also equal to $d$. We have the following result  which is a consequence of results by Gromov and Yomdin \cite{Gromov2,Yomdin}.

\begin{theorem}
The topological entropy of $f$ is equal to $\log d$. 
\end{theorem}

In what follows, we only consider automorphisms $f$ with positive entropy, i.e. with dynamical degree $d>1$. Several examples of classes of such automorphisms can be found in \cite{BK1,Cantat1,DS5,MM, Oguiso1,Uehara}. 
The following proposition describes the action of $f$ on cohomology, see \cite{Cantat1}. 

\begin{proposition} \label{prop_cohomology_action}
There are unique classes $c_+$ and $c_-$ in $H^{1,1}(X,\R)$ such that 
$$c_+\smallsmile \{\omega\}=c_-\smallsmile\{\omega\}=1, \quad f^*(c_+)=dc_+ \quad \mbox{and} \quad f_*(c_-)=dc_-.$$
Moreover, $c_+,c_-$ belong to the boundary of the K\"ahler cone $\Kc$. We have $c_+^2=c_-^2=0$ and $c_+\smallsmile c_-\not=0$. The cup-product $\smallsmile$ is negative definite on the linear space
$$H:=\{c\in H^{1,1}(X,\R),\ c\smallsmile c_+=c\smallsmile c_-=0\}.$$
We have $H^{1,1}(X,\R)=\R c_+\oplus \R c_-\oplus H$. 
The operators $f^*,f_*$ preserve $H$ and act isometrically on $H$ with respect to the cup-product $\smallsmile$. In particular, the actions of $f^*,f_*$ on $H^{1,1}(X,\C)$ are diagonalizable and
$d, d^{-1}$ are their only eigenvalues  with modulus different of $1$.
\end{proposition}
\proof
Recall that $\Kc$ is the cone of classes of K\"ahler forms in $H^{1,1}(X,\R)$. 
Observe that $\overline\Kc$ is a closed strictly convex cone which is invariant under $f^*$ and $f_*$. Since $d$ is the spectral radius of $f^*$ and $f_*$ on $H^{1,1}(X,\R)$, a version of the classical Perron-Frobenius theorem insures the existence of non-zero classes $c_+,c_-$ in $\overline \Kc$ such that $f^*(c_+)=dc_+$ and $f_*(c_-)=dc_-$. 
The fact that $c_\pm$ are in $\overline\Kc\setminus \{0\}$ implies that $c_\pm\smallsmile \{\omega\}\not=0$. 
Therefore, we can normalize these classes so that $c_+\smallsmile \{\omega\}=c_-\smallsmile\{\omega\}=1$. We also deduce from the above invariant relation that $f_*(c_+)=d^{-1}c_+$ and $f^*(c_-)=d^{-1}c_-$. Therefore, $c_+$ and $c_-$ are linearly independent.

Since $f^*=f_*=\id$ on $H^{2,2}(X,\C)$, these operators preserve the cup-product, i.e. we have $f^*(c)\smallsmile f^*(c')=c\smallsmile c'$ for  $c,c'\in H^{1,1}(X,\C)$ and a similar identity holds for $f_*$. We then  deduce from the identity $f^*(c_+^2)=d^2c_+^2$ that $c_+^2=0$. In particular, $c_+$ is on the boundary of $\Kc$. We obtain in the same way that $c_-^2=0$ and that $c_-$ is on the boundary of $\Kc$. It follows from Corollary \ref{cor_HR} that $c_+\smallsmile c_-\not=0$. 

The last property and the fact that $c_+,c_-$ are linearly independent imply that $H$ is a codimension 2 subspace of $H^{1,1}(X,\R)$ which is invariant under $f^*$ and $f_*$.  
We also have  $H^{1,1}(X,\R)=\R c_+\oplus \R c_-\oplus H$. Recall that $H$ is the codimension 2 subspace orthogonal to $c_+$ and $c_-$ with respect to the cup product. 
Theorem \ref{th_Dinh_Nguyen} implies by continuity that the cup-product is semi-negative on $H$. By
Corollary \ref{cor_HR}, 
if $c\in H\setminus\{0\}$, since $c_+^2=0$ and $c_+\smallsmile c=0$, we have $c^2\not=0$. So the cup-product 
is non-degenerate on $H$. Therefore, it is negative definite on $H$. The proposition easily follows.
\endproof

Note that if we replace $\omega$ by another K\"ahler class, then $c_+$ and $c_-$ change by some multiplicative constants. In what follows, for simplicity, we normalize $\omega$ by multiplying it with a constant so that 
$$c_+\smallsmile c_-=1.$$

We have the following result. The simple rigidity of $c_+$ and $c_-$ was proved in \cite{Cantat1} for K3 surfaces and in \cite{DS2} for general surfaces.

\begin{theorem} \label{th_Green_rigid}
The classes $c_+$ and $c_-$ are very rigid. Let $T_+,T_-$ be the unique positive closed $(1,1)$-currents in $c_+$ and $c_-$ respectively. Then $T_+,T_-$ have H\"older continuous local potentials and satisfy
$$f^*(T_+)=dT_+\quad \text{and} \quad f_*(T_-)=dT_-.$$
\end{theorem}
\proof {\bf (except for the H\"older continuity of potentials)} 
We only consider the class $c_+$. The case of $c_-$ can be treated in the same way. Since $c_+$ is in the boundary of $\Kc$, it contains a positive closed current $T_+$. Replacing $T_+$ by a limit of the sequence
$${1\over n} (T_++d^{-1}f^*(T_+)+\cdots + (f^{n-1})^*(T_+))$$
allows to assume that $d^{-1}f^*(T_+)=T_+$. 

Let $S$ be another positive closed current in $c_+$. We prove that $S=T_+$. Write $S_n:=d^n(f^n)_*(S)$. This is also a positive closed current in $c_+$. Let $S_n-T_+=\ddc u_n$ where $u_n$ is a d.s.h. function normalized so that $\int_X u_n\omega^2=0$. So the DSH-norm of $u_n$ is bounded independently of $n$. Let $\varphi$ be a test smooth $(1,1)$-form. Define $\nu:=\ddc\varphi$ and $\nu_n:=(f^n)_*(\nu)$. So $\nu_n$ is a form of maximal degree and defines a measure of bounded mass. Using that $\|f^{-1}\|_{\Cc^1}$ is bounded we obtain that $\|\nu_n\|_\infty\leq M^n$ for a constant $M>0$ large enough.

We have
\begin{eqnarray*}
\langle S-T_+,\varphi\rangle &=& \langle d^{-n}(f^n)^*(S_n)-d^{-n}(f^n)^*(T_+),\varphi\rangle =d^{-n} \langle \ddc (u_n\circ f^n),\varphi\rangle \\
& = & d^{-n} \langle u_n\circ f^n,\ddc \varphi\rangle=d^{-n} \langle u_n\circ f^n,\nu\rangle = d^{-n}\langle\nu_n,u_n\rangle.
\end{eqnarray*} 
By Corollary \ref{cor_dsh_exp}, the last expression is of order $O(nd^{-n})$. Therefore, taking $n\to\infty$, we obtain that $S=T_+$. So $c_+$ is a rigid class.

Consider now a positive $\ddc$-closed current $S$ in $c_+$. We have to check that $S=T_+$. It is enough to prove that $S$ is closed. We only show that $\dbar S=0$ since one can obtain in the same way that $\partial S=0$. Define $S_n$ as above. By Theorem \ref{th_ddc_decom}, we can write for a closed smooth real $(1,1)$-form $\alpha_+$ in $c_+$
$$S_n=\alpha_++\partial\sigma_n+\overline{\partial\sigma_n}$$
where $\sigma_n$ is a $(0,1)$-current  such that $\dbar\sigma_n$ is a $(0,2)$-form with $L^2$-norm bounded independently of $n$. 

We have
$$\dbar S=d^{-n} \dbar (f^n)^*(S_n) =-\partial \big[d^{-n} (f^n)^*(\dbar \sigma_n)\big].$$
On the other hand, we have 
$$\int_X (f^n)^*(\dbar \sigma_n)\wedge (f^n)^*(\partial\overline \sigma_n) = \int_X (f^n)^*(\dbar \sigma_n\wedge \partial\overline \sigma_n)=\int_X \dbar \sigma_n\wedge \partial\overline \sigma_n.$$
We used here the fact that $f$ is an automorphism. So the $L^2$-norm of  $(f^n)^*(\dbar \sigma_n)$ is bounded independently of $n$. Taking $n\to\infty$ gives $\dbar S=0$. We conclude that $c_+$ is very rigid. The proof of the H\"older continuity of potentials of $T_\pm$ will be given below after Theorem \ref{th_explicit_rigidity}.
\endproof
The following result applies in particular to subvarieties $V$ in $X$ and gives an equidistribution property of $f^{-n}(V)$ when $n\to\infty$. 

\begin{theorem} \label{th_equi_Green_surface}
Let $(S_n)$ be a sequence of positive closed $(1,1)$-currents in a given cohomology class $c$. Let $\lambda_c:=c\smallsmile c_-$. Then $d^{-n}(f^n)^*(S_n)\to \lambda_c T_+$ exponentially fast: there is a constant $A>0$ such that 
$$|\langle d^{-n}(f^n)^*(S_n)-\lambda_c T_+,\varphi\rangle| \leq A nd^{-n} \|\varphi\|_{\Cc^2},$$
for every test $(1,1)$-form $\varphi$ of class $\Cc^2$. 
\end{theorem}
\proof
Without loss of generality, we can assume that $\lambda_c=1$. Let $c'$ be a class in $H^{1,1}(X,\C)$ such that $f^*c'=\lambda c'$ with $|\lambda|\leq 1$. Let $\alpha$ be a smooth closed $(1,1)$-form in $c'$. We first show that $d^{-n}(f^n)^*(\alpha)\to 0$ exponentially fast. 

Write $f^*(\alpha)=\lambda \alpha+\ddc u$ with $u$ smooth. We have 
$$d^{-n} (f^n)^*(\alpha)=d^{-n}\lambda^{n} \alpha + \ddc \big[d^{-n}(u\circ f^{n-1}+\cdots+\lambda^{n-1} u)\big].$$
Denote by $u_n$ the function in the brackets. We have $\|u_n\|_\infty=O(nd^{-n})$. It follows that 
$$|\langle d^{-n} (f^n)^*(\alpha),\varphi\rangle| \leq d^{-n}|\langle\alpha,\varphi\rangle|+|\langle \ddc \varphi,u_n\rangle|.$$
It is clear that $d^{-n} (f^n)^*(\alpha)\to 0$ with speed $O(nd^{-n})$. So by Proposition \ref{prop_cohomology_action}, if $\alpha$ is a smooth form whose class is in the hyperplane $H\oplus \R c_-$, then $d^{-n}(f^n)^*(\alpha)\to 0$ with speed $O(nd^{-n})$. 

Since $\lambda_c=1$ we can write $c=c_++c'$ with $c'\in H\oplus \R c_-$. Let $\alpha$ be a smooth closed $(1,1)$-form in $c'$. It is enough to show that $d^{-n} (f^n)^*(S_n-T_+-\alpha)\to 0$ with speed $O(nd^{-n})$. Write $R_n:=S_n-T_+-\alpha=\ddc v_n$ where $v_n$ is a d.s.h. function with DSH-norm bounded independently of $n$. 

Define $\nu:=\ddc \varphi$ and $\nu_n:=(f^n)_*(\nu)$. We obtain as in Theorem \ref{th_Green_rigid} that 
$$\langle d^{-n}(f^n)^*(R_n),\varphi\rangle =d^{-n}\langle \nu_n, v_n\rangle.$$
By Corollary \ref{cor_dsh_exp}, the last expression is of order $O(nd^{-n})$. This completes the proof of the theorem. 
It is possible to weaken the hypothesis on the classes of $S_n$ but one has to take into account the convergence speed of the classes of $d^{-n}(f^n)^*(S_n)$ in $H^{1,1}(X,\R)$. 
\endproof

The above arguments  give the following more quantitative rigidity result for the classes $c_+$ and $c_-$.

\begin{theorem} \label{th_explicit_rigidity}
Let $T_n$ be a sequence of positive closed $(1,1)$-currents. Let $c_n$ denote the cohomology class of $T_n$. Define also $\lambda_n:=\|c_n-c_+\|$ for a fixed norm on $H^{1,1}(X,\R)$. Assume that $c_n\to c_+$, i.e. $\lambda_n\to 0$. Then $T_n\to T_+$ with speed $|\log\lambda_n|\lambda_n^{1/2}$: there is a constant $A>0$ such that 
$$|\langle T_n-T_+,\varphi\rangle| \leq A|\log\lambda_n|\lambda_n^{1/2} \|\varphi\|_{\Cc^2},$$
for every test $(1,1)$-form $\varphi$ of class $\Cc^2$. 
\end{theorem}
\proof
Let $T$ be a positive closed $(1,1)$-current and $c$ its cohomology class. Define $\lambda:=\|c-c_+\|$ and assume that $|\lambda|\ll 1$. It is enough to show for some constant $A>0$ that 
$$|\langle T-T_+,\varphi\rangle| \leq A|\log\lambda|\lambda^{1/2} \|\varphi\|_{\Cc^2}.$$

Multiplying $T$ with a constant close to 1, of order $1+O(\lambda)$,  allows to assume that $c-c_+$ belongs to the hyperplane $H\oplus \R c_-$. 
Let $n$ be the integer part of ${1\over 2}|\log\lambda|/\log d$. Define $S_n:=d^n(f^n)_*(T)=T_++d^n(f^n)_*(T-T_+)$. We have $d^{2n}\simeq \lambda^{-1}$ and  $T=d^{-n}(f^n)^*(S_n)$. Moreover, due to the choice of $n$, the cohomology class $\{S_n\}$ is the sum of $c_+$ and a bounded class in $H\oplus \R c_-$. 
We need here the fact that $\|(f^n)_*\|\lesssim d^n$. 
So the same arguments as in Theorem \ref{th_equi_Green_surface} give the result.
\endproof

\noindent
{\bf End of the proof of Theorem \ref{th_Green_rigid} (H\"older continuity of potentials).} Let $\alpha_+$ be a smooth real $(1,1)$-form in $c_+$. Since we can write it as a difference of two K\"ahler forms, by Theorem \ref{th_equi_Green_surface}, $d^{-n}(f^n)^*(\alpha_+)$ converge to a constant times $T_+$. Since $\alpha_+$ is in $c_+$, this constant should be 1. So we have $d^{-n}(f^n)^*(\alpha_+)\to T_+$. 

We can write for some smooth function $v^+$ that $d^{-1}f^*(\alpha_+)=\alpha_++\ddc v^+$. A simple induction on $n$ gives
$$d^{-n} (f^n)^*(\alpha_+)=\alpha_+ + \ddc (v^++d^{-1}v^+\circ f+\cdots + d^{1-n} v^+\circ f^{n-1}).$$ 
Taking $n\to\infty$ gives $T_+=\alpha_+ +\ddc v^+_\infty$ with $v^+_\infty:=\sum_{n\geq 0} d^{-n}v^+\circ f^n$. Lemma \ref{lemma_Holder} implies that $v^+_\infty$ is H\"older continuous. So $T_+$ has H\"older continuous local potentials. This completes the proof of the theorem. 
\hfill $\square$

\begin{corollary}
Let $S$ be a positive $(1,1)$-current on $X$ such that $\ddc S$ is of order $0$. Then we have $d^{-n} (f^n)^*(S)\to c_ST_+$, where $c_S:=\langle T_-,S\rangle$. 
\end{corollary}
\proof
Replacing $S$ with a suitable combination of $S$ and $\omega$ allows us to  assume that $c_S=1$. 
Define $S_n:=d^{-n}(f^n)^*(S)$. We first show that $\|S_n\|\to 1$. As a consequence of Theorem \ref{th_equi_Green_surface} applied to $f^{-1}$ instead of $f$, we have $d^{-n}(f^n)_*(\omega)\to T_-$.  Let $\alpha_-$ be a smooth real closed $(1,1)$-form in $c_-$. As in the proofs of Theorems \ref{th_Green_rigid} and \ref{th_equi_Green_surface}, we can write 
$d^{-n}(f^n)_*(\omega)-\alpha_-$ as a sum  $\beta_n+\ddc v_n^-$, where $\beta_n$ are smooth forms such that $\|\beta_n\|_\infty\to 0$ and $v_n^-$ are smooth functions converging uniformly to a continuous  function $v^-_\infty$ such that $T_-=\alpha_-+\ddc v^-_\infty$.  

We have 
\begin{eqnarray*}
\|S_n\| & = & \langle \omega, d^{-n}(f^n)^*(S)\rangle =\langle d^{-n}(f^n)_*(\omega),S\rangle \\
& = & \langle \alpha_-,S\rangle +\langle\beta_n,S\rangle +\langle \ddc v^-_n,S\rangle = \langle \alpha_-,S\rangle +\langle\beta_n,S\rangle +\langle \ddc S,  v^-_n\rangle.
\end{eqnarray*}
Clearly, $\|S_n\|$ converge to $\langle T_-,S\rangle=1$.

By hypothesis, $\ddc S_n$ is a measure. Since $f$ is an automorphism, we have 
$$\|\ddc S_n\|=d^{-n}\|(f^n)^*(\ddc S)\|=d^{-n}\|\ddc S\|\to 0.$$
So the sequence of currents $S_n$ is relatively compact and the family $\Fc$ of cluster values contains only positive $\ddc$-closed currents of mass 1. By definition, this family is compact and invariant under $d^{-1}f^*$. 

The set $\Gc$ of classes $\{T\}$ with $T\in\Fc$ is then compact, disjoint from 0 and invariant under $d^{-1}f^*$. The description of $f^*$ in Proposition \ref{prop_cohomology_action} implies that $\Gc$ is contained in the half-line generated by $c_+$. So by Theorem \ref{th_Green_rigid},  $T$ is proportional to $T_+$. Since its mass is equal to 1, we necessarily have $T=T_+$. This completes the proof of the corollary.
\endproof

\begin{remark} \rm
We can prove using similar arguments that $d^{-n}(f^n)^*(S_n)\to T_+$ if $S_n$ are positive $(1,1)$-currents such that $\|S_n\|=o(d^n/n)$, $\|\ddc S_n\|=o(d^n)$ and $\langle T_-,S_n\rangle\to 1$. 
\end{remark}

The following result can be applied to stable manifolds.

\begin{corollary} \label{cor_C_T}
Let $\phi:\C\to X$ be a holomorphic map such that $\phi(\C)$ is not contained in a proper subvariety of $X$. Assume that $\phi^*(T_+)=0$. Then $T_+$ is the only Nevanlinna current  associated with $\phi$. 
In particular, the currents $\tau_r$ defined in Definition \ref{def_Ne_current} converge to $T_+$.
\end{corollary}
\proof
Recall that if $S$ is a positive closed $(1,1)$-current such that $S=\alpha+\ddc u$ with $\alpha$ smooth and $u$ continuous quasi-p.s.h., we define $\phi^*(S):=\phi^*(\alpha)+\ddc (u\circ\phi)$. It is easy to check that this is a positive measure which does not depend on the choice of $\alpha$ and $u$. 
If we write locally $S=\ddc v$ with $v$ continuous p.s.h., then $\phi^*(S)=\ddc(v\circ\phi)$ on the domain of definition of $v\circ\phi$. 
The hypothesis $\phi^*(T)=0$ means that local potentials of $T$ are harmonic on the image of $\phi$.

Recall also that  $\tau_r$  is a positive current such that $\|\ddc \tau_r\|$ bounded independently of $r$. We have 
$$\langle T_+,\tau_r\rangle =\Big\langle T_+, {1\over T(r)} \int_0^r{dt\over t} \phi_*[\D_t]\Big\rangle= {1\over T(r)}\int_0^r{dt\over t} \int_{\D_r} \phi^*(T_+)=0.$$
Let $T$ be a Nevanlinna current associated with $\phi$. Recall that $T$ is positive $\ddc$-closed of mass $1$.
We deduce from the above identities and properties of intersection number introduced in Section \ref{section_Kahler} that  $\langle T_+,T\rangle=0$ and hence $\{T\}\smallsmile c_+=0$. 

A version of McQuillan's theorem \cite{BS,MQ} says that $\{T\}$ is nef, i.e. is in $\overline\Kc$. Therefore, we have $\{T\}^2\geq 0$. Corollary \ref{cor_HR} implies that $\{T\}$ is proportional to $c_+$. Since $\|T\|=1$, by Theorem \ref{th_Green_rigid}, we necessarily have $T=T_+$. The proof is the same for  Nevanlinna currents associated with Riemann surfaces which are not necessarily equal to $\C$. 
\endproof

As in Corollary \ref{cor_Fatou_Henon}, we deduce from the last result  that 
if a Fatou component $U$ of $f$ is  not hyperbolically embedded in $X$ then its boundary contains the support of $T_+$. Indeed, in this case, there is a sequence of holomorphic discs $\phi_n:D(0,n)\to U$ converging to a non-constant map $\phi:\C\to \overline U$. Since $\phi_n^*(T)=0$, we also have $\phi^*(T)=0$. 
This is the case for basins of attracting points, see \cite{MM} for examples and \cite{Moncet} for a related result. 
The following proposition characterizes  curves for which Corollary \ref{cor_C_T} applies.

\begin{proposition}
Let $\phi:\Sigma\to X$ be a holomorphic map on a Riemann surface $\Sigma$. Then $\phi^*(T_+)=0$ if and only if for any compact set $K\subset\Sigma$ the area of  $f^n(\phi(K ))$ counted with multiplicity is equal to $O(n)$ as $n\to\infty$. In particular, we have $\phi^*(T_+)=0$ when the sequence $(f^n\circ\phi)_{n\geq 0}$ is locally equicontinuous on  $\Sigma$. 
\end{proposition}
\proof
When the sequence $(f^n\circ\phi)_{n\geq 0}$ is locally equicontinuous on  $\Sigma$, the area of $f^n(\phi(K ))$ is bounded independently of $n$. So the second assertion is a consequence of the first one.

We have 
$$\phi^*(T_+)=\lim_{n\to\infty} \phi^*(d^{-n}(f^n)^*(\omega))=\lim_{n\to\infty} d^{-n}(f^n\circ\phi )^*(\omega).$$
So if the area of  $f^n(\phi(K ))$ is equal to $o(d^n)$ for every $K$, then $\phi^*(T_+)=0$.

Assume now that $\phi^*(T_+)=0$. We have to show that the area of   $f^n(\phi(K ))$ is equal to $O(n)$.
Choose a smooth function $0\leq \chi\leq 1$ with compact support in $\Sigma$ and equal to 1 on $K$. With the above notation, we can write 
$$\omega=T_++\alpha+\ddc u,$$
where $u$ is a continuous d.s.h. function and $\alpha$ is a smooth $(1,1)$-form in a class of $H\oplus \R c_-$. 
We have
$$(f^n)^*(\omega)=d^nT_++(f^n)^*(\alpha)+\ddc (u\circ f^n).$$
We can write as in Theorem \ref{th_equi_Green_surface} 
$$(f^n)^*(\alpha)=\alpha_n+\ddc v_n$$
with $\alpha_n$ bounded uniformly on $n$ and $\|v_n\|_\infty=O(n)$.

So the area of $f^n(\phi(K))$ is bounded by
$$\int_\Sigma \chi \phi^*(f^n)^*(\omega)=d^n\langle\phi^*(T_+),\chi\rangle +\langle \phi^*(\alpha_n),\chi\rangle +\langle \phi_*(\ddc \chi), v_n+u\circ f^n\rangle.$$
Since the first term vanishes, it is now clear that the area of $f^n(\phi(K))$ is equal to $O(n)$ as $n\to\infty$. 
\endproof

Finally, we have the following result.

\begin{proposition}
Assume that $c_++c_-$ is a K\"ahler class. Then, the  support of $T_+$ is equal to the forward Julia set. 
\end{proposition}
\proof
Let $F'$ be the complement of $\supp(T_+)$. We can assume that $\omega$ is a K\"ahler form in $c_++c_-$. 
We know that $T_+$ is the limit of $d^{-n} (f^n)^*(\omega)$. On the Fatou set $F$, the forms $(f^n)^*(\omega)$ are locally bounded uniformly on $n$. Therefore, we have $F\subset F'$. 

We can write 
$$\omega=T_++T_-+\ddc u$$
where $u$ is a continuous d.s.h. function. We have on $F'$
$$(f^n)^*(\omega)=d^{-n}T_-+\ddc (u\circ f^n).$$ 
In particular, $(f^n)^*(\omega)$ admits in $F'$ local potentials which are bounded uniformly on $n$. Therefore, this family of currents is relatively compact and their cluster values have locally bounded potentials on $F'$. 
The lemma below implies that $(f^n)$ is locally equicontinuous on $F'$. This completes the proof of the proposition. 
Note that the volume of the graph of $f^n$ over $F$, $n\geq 0$,  is locally bounded independently of $n$. 
\endproof

The following lemma is essentially obtained in \cite{Dinh2}, see also \cite{DS7}. 

\begin{lemma}
Let $(h_n)$ be a sequence of holomorphic maps from a complex manifold $U$ to a fixed compact subset $K$ of a K\"ahler manifold $V$. Let $\omega$ be a K\"ahler form on $V$. Assume that the family of positive closed $(1,1)$-currents $h_n^*(\omega)$ is relatively compact and for any cluster value $S$ of this sequence the Lelong number of $S$ at every point is bounded by a fixed constant $c(V,K,\omega)$ small enough. Then the family $(h_n)$ is locally equicontinuous. In particular, if $L$ is a compact subset of $U$ then the volume of $h_n(L)$ is bounded independently of $n$.  
\end{lemma}

We refer to \cite{Demailly3} for the notion of Lelong number. We only need here the fact that the Lelong number of a current with bounded local potentials always vanishes. Note that without the hypothesis on the Lelong number, we can extract from $(h_n)$ a subsequence which converges locally uniformly outside the analytic set of points where the Lelong number of $S$ is larger than $c(V,K,\omega)$. This is a higher dimensional version of a famous lemma due to Gromov which is valid for maps defined on a Riemann surface.

\begin{remark} \rm
Of course the above results can be applied to $f^{-1}$ and to $T_-$. The intersection $\mu:=T_+\wedge T_-$ defines an invariant probability measure. This measure turns out to be the unique measure of maximal entropy $\log d$. It is exponentially mixing and when $X$ is projective saddle periodic points are equidistributed with respect to $\mu$. When $X$ is a projective surface, the techniques developed for H\'enon type maps can be applied without difficulty. In general, new ideas and even completely new tools are needed. Several properties including the statistical ones still hold for large classes of horizontal-like maps. We refer to \cite{Cantat1, Cantat2, deThelin1, dD, Dinh1, DNS2,DS3, DS10, DS9, Du} for these developments.
\end{remark}

\section{Dynamics in higher dimension} \label{section_higher_dim}

In this section, we briefly discuss similar situations in higher dimension: polynomial automorphisms of $\C^k$ and holomorphic automorphims of compact K\"ahler manifolds.

Let $f$ be a polynomial automorphism of $\C^k$. We still denote by $f$ its extension as a birational map of $\P^k$. Let $I_+,I_-$ denote the indeterminacy sets of $f$ and $f^{-1}$ respectively. 
They are analytic sets strictly contained in the hyperplane at infinity. We assume that $I_+$ and $I_-$ are non-empty; otherwise, $f$ is an automorphism of $\P^k$ and its dynamics is easy to understand.
The following notion was introduced by the second author in \cite{Sibony}.

\begin{definition} \rm
We say that $f$ is {\it regular or of H\'enon type} if $I_+\cap I_-=\varnothing$.
\end{definition}

It is remarkable that the later condition is quite easy to check while it should be difficult to develop a theory  for all automorphisms of $\C^k$ with $k\geq 3$, see e.g. \cite{SU}.
Moreover, the family of regular maps is very rich. In dimension 2, we have seen that any polynomial automorphism dynamically interesting is conjugated to a regular one. 

From now on, assume that $f$ is regular.
Denote by $d_+,d_-$ the algebraic degrees of $f$ and $f^{-1}$ respectively. We recall here some elementary properties of $f$ and refer to \cite{Sibony} for details.

\begin{proposition}
There is an integer $1\leq p\leq k-1$ such that $\dim I_+=k-p-1$, $\dim I_-=p-1$ and $d_+^p=d_-^{k-p}$. We have 
$$f(\{w_0=0\}\setminus I_+)=f(I_-)=I_- \quad \text{and}\quad  f^{-1}(\{w_0=0\}\setminus I_-)=f^{-1}(I_+)=I_+.$$  
In particular, $f^n$ is a regular automorphism for every $n\geq 1$ and the algebraic degrees of $f^n, f^{-n}$ are $d_+^n,d_-^n$ respectively. 
\end{proposition}

\begin{proposition}
The set $I_-$ (resp. $I_+$) is attracting for $f$ (resp. for $f^{-1}$). Denote by $\widetilde U_+$ (resp. $\widetilde U_-$) its basin. Define also $U_\pm:=\widetilde U_\pm\cap \C^k$ and $K_\pm:=\C^k\setminus U_\pm$.  Then $K_+$ (resp. $K_-$) is the set of points $z\in\C^k$ such that the orbit $(f^n(z))_{n\geq 0}$  (resp. $(f^{-n}(z))_{n\geq 0}$) is bounded in $\C^k$. Moreover, we have $\overline K_\pm=K_\pm\cup I_\pm$. 
\end{proposition}

The Green function for $f$ is defined as in the case of dimension 2. Put 
$$G^+_n(z):=d_+^{-n} \log^+\|f^n(z)\|.$$
The following theorem is obtained in the same way as in the  dimension 2 case.

\begin{theorem}
The sequence $G_n^+$ converges locally uniformly on $\C^k$ to a H\"older continuous p.s.h. function $G^+$ such that $G^+(z)-\log^+\|z\|$ extends to a H\"older continuous function on $\P^k\setminus I_+$. Moreover, we have $G^+\circ f=dG^+$. 
\end{theorem}

As a consequence of the above theorem, we obtain the following result.

\begin{corollary}
The current $T_+:=\ddc G^+$ extends by zero to a positive closed $(1,1)$-current of mass $1$ on $\P^k$. Moreover, for $1\leq q\leq p$, the power $T_+^q$ is well-defined and is a positive closed $(q,q)$-current of mass $1$ on $\P^k$. We also have $f^*(T_+^q)=d_+^q T_+^q$ for $1\leq q\leq p$ and $\supp(T_+^p)\subset J_+:=\partial\overline K_+$. 
\end{corollary}

\begin{definition} \rm
We call $G^+$ {\it the Green function}  of $f$ and $T_+^q$ {\it the Green $(q,q)$-current} of $f$.
\end{definition}

The following results were obtained by the authors in \cite{DS6}. The proofs are much more delicate than in the dimension 2 case. They are based on a theory of super-potentials which allows to deal with positive closed $(p,p)$-currents, $p>1$.

\begin{theorem}
The set $\overline K_+$ is $p$-rigid and $T_+^p$ is the unique positive closed $(p,p)$-current of mass $1$ supported on $\overline K_+$.
\end{theorem}

We have the following strong equidistribution property. It can be applied to currents of integration on analytic sets. For the case of bidegree $(1,1)$ see \cite{DS6, Taflin}. 

\begin{theorem}
Let $U$ be a neighbourhood of $I_-$. Let $S_n$ be a sequence of positive closed $(p,p)$-currents of mass $1$ with support in $\P^k\setminus U$. Then $d_+^{-pn} (f^n)^*(S_n)\to T_+^p$ exponentially fast.
\end{theorem}

We however don't know if $\overline K_+$ is very $p$-rigid except for $p=1$.

\begin{theorem}
Assume that $p=1$. Then $\overline K_+$ is very rigid.
\end{theorem} 
\proof
The proof follows the one of Theorem \ref{th_K_rigid}. We use the same notation. The only different point is the estimation of the $L^2$-norm of $(f^n)^*(\dbar\sigma_n)$.  We need to show that this norm is equal to $o(d_+^n)$. 
Observe that since $S_n$ is supported on $\overline K_+$, the forms $\dbar\sigma_n$ and $\partial\overline\sigma_n$ are smooth near $I_-$. 
We have
$$\int_{\P^k} (f^n)^*(\dbar\sigma_n) \wedge (f^n)^*(\partial\overline\sigma_n) \wedge\omega_\FS^{k-2} =
\int_{\P^k} \dbar\sigma_n\wedge\partial\overline\sigma_n \wedge (f^n)_*(\omega_\FS^{k-2}).$$

On the other hand, $(f^n)_*(\omega_\FS^{k-2})$ is positive closed and smooth outside $I_-$. Its mass is equal to $d_-^{(k-2)n}$. Therefore, by Theorem \ref{th_ddc_pot}, the last integral is of order $O(d_-^{(k-2)n})=o(d_+^{n})$ since $d_+=d_-^{k-1}$. The rest of the proof is the same as in Theorem \ref{th_K_rigid}.
\endproof

\begin{remark} \rm
The automorphism $f^{-1}$ is also regular of algebraic degree $d_-$. We can construct as above the Green function $G^-$ and the Green $(q,q)$-currents $T_-^q$ for $f^{-1}$ with $1\leq q\leq k-p$. The Green current $T_-^{k-p}$ is the unique positive closed $(k-p,k-p)$-current of mass 1 supported on the rigid set $\overline K_-$. The measure $\mu:=T_+^p\wedge T_-^{k-p}$ is the unique invariant measure of maximal entropy $p\log d_+$. It is exponentially mixing and saddle periodic point are equidistributed with respect to $\mu$. We refer to \cite{deThelin2, Dinh1, DS6, DS9} for details. Note that for simplicity the exponential mixing was only given in \cite{Dinh1} under the hypothesis $k=2p$ but there is no difficulty to extend it to the general case, see also \cite{DS10} for the techniques which give the best estimates using this approach. The main new tools which allow to deal with dynamics in higher dimension were developed by the authors in \cite{DS1,DS6,DS2}.
\end{remark}

For the rest of this section, we discuss the case of automorphisms of compact K\"ahler manifolds. Let $f:X\to X$ be an automorphism on a compact K\"ahler manifold $(X,\omega)$ of dimension $k$. 

\begin{definition} \rm
We call {\it dynamical degree of order $q$} of $f$ the spectral radius of $f^*$ on $H^{q,q}(X,\R)$. 
\end{definition}

It is not difficult to see that $d_0=d_k=1$.
The following result is a consequence of Theorem \ref{th_Dinh_Nguyen}, results by Gromov and Yomdin \cite{Gromov2,Yomdin} and some observations from \cite{DS5}.

\begin{theorem}
The dynamical degrees of $f$ are log-concave in $q$, i.e. $d_q^2\geq d_{q-1}d_{q+1}$ for $1\leq q\leq k-1$. In particular, there are two integers $0\leq  p,p'\leq k$ such that $d_0<\cdots <d_p=\cdots=d_{p'}>\cdots>d_k$. Moreover, the topological entropy of $f$ is equal to $\log d_p$. In particular, $f$ has positive entropy if and only if $d_1>1$. In this later case, the entropy of $f$ is bounded below by a positive constant depending only on the second Betti number of $X$. 
\end{theorem}

The relation with the second Betti number is due to the following fact. The topological entropy is bounded from below by the logarithm of the spectral radius of $f^*$ on $H^2(X,\C)$. The last number is the largest modulus of the roots of a polynomial with integer coefficients whose degree is the second Betti number. If one of the coefficients of this polynomial is too large, it is clear that there is a big root. Otherwise, we have a finite number of polynomials to consider and the result follows easily. One can also consider the other Betti number with even index.

\smallskip

The following rigidity theorem was obtained by the authors in \cite{DS2} using the theory of super-potentials.

\begin{theorem}
Let $E\subset H^{q,q}(X,\R)$,  with $1\leq q\leq p$, be a linear subspace invariant under $f^*$. Assume that all complex eigenvalues of $f^*_{|E}$ have modulus strictly larger than $d_{q-1}$. If a class $c$ in $E$ contains a non-zero positive closed $(q,q)$-current, then $c$ is rigid. 
\end{theorem}

Note that $d_q$ is an eigenvalue of $f^*$ on $H^{q,q}(X,\R)$. So we can construct positive closed $(q,q)$-currents $T_+$ in some rigid cohomology classes  such that $f^*(T_+)=d_qT_+$. We call them {\it Green $(q,q)$-currents}. They have H\"older continuous super-potentials, see  \cite{DS2} for details.

Equidistribution results can be deduced from the last theorem by observing that if $S_n$ are positive closed currents such that $\{S_n\}$ converge to the rigid class $\{T\}$ of a positive closed current $T$, then $S_n$ converge to $T$.  If $T$ is a Green current, the speed of convergence of $S_n$ can be bounded in term of the speed of convergence of $\{S_n\}$. However, the action of $f^*$ on cohomology is far from being well-understood. The presence of Jordan blocks may induce slow convergence. Besides of this difficulty, the control of the convergence speed is satisfactory. The following result can be deduced from our study in \cite{DS2}. The estimate can be improved using eigenvalues of $f^*$ on $H^{q-1,q-1}(X,\R)$ and on $H^{q,q}(X,\R)$. 

\begin{theorem}
Let $T_+$ be a Green $(q,q)$-current of $f$ and $c_+$ its cohomology class. 
Let $S$ be a positive closed $(q,q)$-current  and $c$ its cohomology class such that $\lambda:=\|c-c_+\|$ is small enough. Then, there
are constants $A>0$ and $\alpha>0$ independent of $S$ such that 
$$|\langle S-T_+,\varphi\rangle| \leq A \lambda^\alpha \|\varphi\|_{\Cc^2},$$
for every test $(k-q,k-q)$-form $\varphi$ of class $\Cc^2$. 
\end{theorem}

We don't know if the cohomology class of a Green current is always very rigid except for $q=1$. The following result is obtained as in the case of surfaces and the case of regular automorphisms of $\C^k$.

\begin{theorem}
Assume that $d_2<d_1^2$. Then $f$ admits only one Green current of bidegree $(1,1)$ up to a multiplicative constant. Moreover, its cohomology class is very rigid.
\end{theorem}

Under appropriate conditions on the action of $f$ on cohomology, we can construct for $f$ an invariant probability measure $\mu$ which turns out to be the unique measure of maximal entropy $\log d_p$. This measure is exponentially mixing. We refer to \cite{dD,DS1,DS2,DS10} for details. In a forthcoming work, we will study the equidistribution of  saddle periodic points  with respect to $\mu$. As in the polynomial case, the new tools allowing to treat the higher dimension case were developed in  \cite{DS2,DS6,DS9}.

\noindent
T.-C. Dinh, 
Department of Mathematics, National University 
of Singapore, 10 Lower Kent Ridge Road, Singapore 119076.
 {\tt matdtc@nus.edu.sg}

\

\noindent
N. Sibony,
Universit{\'e} Paris-Sud, Math{\'e}matique - B{\^a}timent 425, 91405
Orsay, France. {\tt nessim.sibony@math.u-psud.fr} 
 \end{document}